\documentclass{amsart}
\usepackage{a4wide}
\usepackage{amsfonts,amsmath,amsthm,amssymb}
\usepackage[linesnumbered,ruled,noend]{algorithm2e}
\usepackage{comment}
\usepackage{float}
\usepackage{latexsym}
\usepackage{graphicx}
\usepackage{multicol}
\usepackage{color}
\usepackage{xcolor}
\usepackage{pifont}
\usepackage{enumerate}
\usepackage{lmodern}
\usepackage{makecell}
\usepackage{adjustbox}
\usepackage{hyperref}
\usepackage{setspace}
\usepackage[all]{xy}
\usepackage{pb-diagram,pb-xy}
\usepackage{pdflscape}
\usepackage{chngpage}
\usepackage{breakurl}
\usepackage{adjustbox,lipsum}
\usepackage{graphicx}
\usepackage[justification=centering]{caption}
\usepackage{float}
\usepackage{dsfont}
\usepackage{ctable}
\usepackage{csquotes}
\usepackage{tabu}

\makeatletter
\def\BState{\State\hskip-\ALG@thistlm}
\makeatother

\theoremstyle{plain}
\newtheorem{theorem}{Theorem}[section]

\newtheorem{conjecture}[theorem]{Conjecture}

\theoremstyle{definition} 

\newtheorem{example}{Example}
\newtheorem*{example*}{Example}

\newcommand{\bbZ}{\mathbb{Z}}
\newcommand{\bbQ}{\mathbb{Q}}

\newcommand{\bbN}{\mathbb{N}}

\newcommand{\Q}{\mathbb{Q}}
\newcommand{\Z}{\mathbb{Z}}
\newcommand{\C}{\mathbb{C}}

\begin{document}

\title[]{Power values of power sums: a survey}

\author{Nirvana Coppola}
\address{Vrije Universiteit Amsterdam, de Boelelaan 1111, Room 9A94, 1081 HV, Amsterdam, The Netherlands}
\email{nirvanac93@gmail.com}

\author{Mar Curc\'{o}-Iranzo}
\address{Hans-Freudental Gebouw, Utrecht University, Budapestlaan 6, Room 5.03, 3584 CD Utrecht, The Netherlands}
\email{m.curcoiranzo@uu.nl}

\author{Maleeha Khawaja}
\address{School of Mathematics and Statistics, University of Sheffield, Hounsfield Road, Sheffield S3 7RH, United Kingdom}
\email{mkhawaja2@sheffield.ac.uk}

\author{Vandita Patel}
\address{School of Mathematics, University of Manchester, Oxford Road, Manchester M13 9PL, United Kingdom}
\email{vandita.patel@manchester.ac.uk}

\author{\"{O}zge \"{U}lkem}
\address{Galatasaray University, \c{C}{\i}ra\u{g}an Cd. No:36, Istanbul, Turkey}
\email{ozgeulkem@gmail.com}

\date{\today}


\begin{abstract}
Research on power values of power sums has gained much attention of late, partially due to the explosion of  refinements in multiple advanced tools in (computational) Number Theory
in recent years.
In this survey, we present the key tools and techniques employed thus far
in the (explicit) resolution of 
Diophantine problems, as well as 
an overview of existing results.
We also state some open problems that naturally arise in the process.

\end{abstract}
\maketitle

\section{Introduction}

Euler, in 1770, noted the remarkable identity
\[
3^3 + 4^3 + 5^3 = 6^3
\]
 in his famous \textit{Vollst\"{a}ndige {A}nleitung zur {A}lgebra} \cite[art. 249]{EulerAlg}. He further asked if there were other instances of sums of consecutive cubes that are perfect powers. Dickson's encyclopaedic \emph{History of the Theory of Numbers} \cite[pp. 582 -- 588]{Dickson}
 extensively surveys this problem citing contributions from many 19th century renowned mathematicians such as Cunningham, Catalan, Genocchi and Lucas. 
 
\medskip

Consecutive cubes that sum to a 
square have the well-known parametric family of solutions
\[
\sum_{i=1}^{d} i^3  = \left(\frac{d(d+1)}{2}\right)^2.
\]
The question of existence of solutions beyond this family 
poses a more interesting and challenging problem. Further studies began with Catalan,  Cantor and Richaud, 
who produced incomplete results by generating solutions via Pell's equation \cite[pp. 586]{Dickson}.

\medskip

Lucas proceeded to boldly make the following three claims (\cite[pg. 92]{Lucas1961}). 

\begin{itemize}
    \item[(a)]  For $d=2$: ``the only square that is the sum of two consecutive positive cubes is $3^{2} = 1^{3} + 2^{3}$''. 
    \item[(b)]  For $d=3$: ``the sum of 3 consecutive cubes is never a square, except for the solution $1^{3}+2^{3}+3^{3}=6^{2}$''. 
    \item[(c)] For $d=5$: ``the sum of 5 consecutive cubes is never a square, except for in the cases where the middle cube is either 2, 3, 27, 98 and 120.'' 
\end{itemize}
Claim (b) turned out to be faulty: Cassels \cite{Cassels} and Uchiyama \cite{Uchiyama} independently prove that the complete list of solutions are actually
\[
    0 = (-1)^{3} + 0^3 + 1^3, \quad 
    3^2 = 0^3 + 1^3 + 2^3, \quad 1^{3}+2^{3}+3^{3}=6^{2}, \quad 
    204^2 = 23^3 + 24^3 + 25^3.
\]
On the other hand, Stroeker \cite{Stroeker1995} proves  that  claims (a) and (c) are indeed correct!  In fact, Stroeker studies the problem of squares that are sums of at most $50$ consecutive cubes using a (by now) standard method based on linear
forms in elliptic logarithms. Nowadays, the problem of finding squares that are sums of $d$ consecutive cubes for a fixed integer $d$ is more commonly viewed as a problem of determining all integral points on certain genus 1 curves \cite{Bennett2016}.

\medskip

Consecutive cubes that sum to a cube possess the lesser-known parametric family of solutions
\[
 \left(\frac{m^{5} + m^3 - 2m}{6}\right)^3 = \sum_{i=1}^{m^3} \left(\frac{m^4 - 3m^3 - 2m^2 - 2}{6} + i\right)^3
\]
where $m \equiv 2 \ \mathrm{ or } \ 4 \pmod{6}$ to guarantee that the cubes are integral. The construction was given in 1829 by Pagliani, whose apparent motivation was to answer a challenge posed presumably by the journal editor Gergonne: to find 1000 consecutive cubes whose sum is a perfect cube. 
Notable mathematicians have contributed to providing parametric families of solutions, including Lenhart and Matthiessen, or proving that in certain instances there are no solutions such as Lebesgue, Evans and Escott. Once more, finding solutions 
outside of Pagliani's family is the more compelling question. 

\medskip

The literature surrounding the general problem of expressing perfect powers as a sum of (a small number of) consecutive cubes (or higher powers) is far more sparse, most likely due to the fact that modern 
 tools and techniques are necessary, even in the case of partial resolution of problems. Since we aim to find all integer solutions to certain exponential Diophantine equations, we can naturally view this as a subset of all of the rational solutions of the same equation. 
 In other words, we are interested in the set of rational points of the algebraic curve $C$ defined by $f(x, y) \in \Z[x, y]$ obtained from the given equation. If $C$ is smooth, Faltings' theorem \cite{Faltings} (which awarded him the Fields Medal in 1986) determines the cardinality of the set of rational points according to the 
 genus of $C$.

\begin{itemize}
    \item If the genus of $C$ is 0, then $C$ is isomorphic to $\mathbb{P}^{1}$ and  $C(\Q)$ is either empty or has infinitely many points parameterised by $\mathbb{P}^{1}(\Q)$.
    \item If the genus of $C$ is 1, and $C$ has a rational point, then $C$ is an elliptic curve and $C(\Q)$ is isomorphic to a finitely generated abelian group $\Z^{r}\oplus T$, where $T$ is finite. 
    \item If the genus of $C$ is greater than or equal to 2, then the set $C(\Q)$ is finite.
\end{itemize}
In particular, note that fixing the genus (equivalent to fixing certain exponents) of an exponential Diophantine equation enables us to conclude, via Faltings' theorem, that there are finitely many integral solutions. However, we encounter two significant drawbacks: firstly,  Faltings' proof is ineffective, thus we are unable to explicitly resolve problems such as those posed by Euler and Gergonne, and secondly, exponential Diophantine equations usually do not have a fixed genus and so we are unable to conclude finiteness results. An effective version of Faltings' proof would indeed be momentous, and certainly pave the way for significant breakthroughs
in the explicit resolution of Diophantine equations!

\medskip

 In this survey, we aim to give a comprehensive
overview of results towards the following families of exponential Diophantine equations:

\begin{equation} \label{eq:table1} 
x^{k}+(x+r)^{k}+ \cdots + (x+(d-1)r)^{k} =y^{n}, 
\qquad d, k, n, r, x, y, \in \Z, \quad 
\end{equation}
with parameters $(k, d, r, n)$.

\begin{equation} \label{eq:table2} s(1^{k}+2^{k}+ \cdots + x^{k}) + r(x) =y^{n}, 
\qquad k, n, x, y, s \in \Z,  r(x) \in \Z[x], \end{equation}
with parameters $(k, n, s, r(x))$.

\medskip

 We begin by highlighting some of the key pioneering advances in number theory that have enabled explicit resolutions of difficult multi-parametric exponential Diophantine equations. 
  One such strategy used to reach a full explicit resolution is to first bound the exponent(s).
  This reduces the number of unknowns in the Diophantine equation.
  The next step is to eliminate equations by ruling out the existence of integral solutions for certain combinations of variables.
  The remaining handful of equations are then resolved using ad hoc techniques.

\subsection*{Bounding the exponent} We present techniques used to bound exponents in certain exponential Diophantine equations.

\subsection{Diophantine approximation via linear forms in logarithms}\label{subsec:linlogs}
Let $\alpha_{1}, \alpha_{2},\dots, \alpha_{n}$ and 
$\beta_{1}, \beta_{2},\dots, \beta_{n}$ be algebraic numbers. A linear form in logarithms is an expression of  the form
\begin{equation}
    \label{lfl}
        \beta_{1}\log\alpha_{1}+\cdots+\beta_{n}\log\alpha_{n}.
\end{equation}
The motivation behind (see \cite[Section 12.1.2]{Cohen}) obtaining a bound on \eqref{lfl} follows from the following idea.
After successfully associating a linear form in logarithms to a putative solution to a given Diophantine equation $D$, the successive step is to then obtain a bound $B_{D}$ for the $\alpha_{i}, \beta_{i}$. 
This provides an upper bound for the solutions to $D$. If $B_{D}$ is small enough then one can attempt to determine all solutions to $D$. 
We note that for Diophantine applications, it is essential that $\beta_{i}\in\bbZ$. Moreover, if the $\alpha_{i}$ are known then one should be able to find a smaller bound on $\beta_{i}$ than $B_{D}$. 

\medskip

The first effective bound was obtained by Baker \cite{Baker66} in 1966, awarding him the Fields medal. 
We note that this is the most general result as it holds for arbitrary $n$, see \cite[pg. 257]{baker90}. 
\begin{theorem}[Baker] \label{thm:Baker}
Let $K$ be a number field. 
Let $\alpha_1, \ldots, \alpha_n \in K^*$ and $\beta_1, \ldots, \beta_n \in K$. For any constant $\nu$, define 
\begin{align*}
\tau(\nu) :&= \tau(\nu;\alpha_1, \ldots, \alpha_n ;\beta_1, \ldots, \beta_n)\\
&=h([1,\beta_1, \ldots, \beta_n])h([1, \alpha_1, \ldots, \alpha_n ])^{\nu}.
\end{align*}
where $h$ is a logarithmic height function. Fix an embedding $K \subset \C$ and let $|\cdot|$ be the corresponding absolute value. Assume that 
\[
\beta_1 \log \alpha_1 + \cdots  +\beta_n \log \alpha_n \neq 0.
\]

Then there are effectively computable constants $C, \nu >0$ which depend only on $n$ and $[K:\Q]$ such that
\[
|\beta_1 \log \alpha_1 + \cdots  +\beta_n \log \alpha_n| > C^{-\tau(\nu)}.
\]
\end{theorem}

\medskip

Baker's theory, albeit effective, usually produces astronomical bounds. Substantial work has been undertaken to reduce such bounds, 
with the majority of the focus 
over the last century residing  within the realm of linear forms in 
\emph{two} logarithms.
One such result is due to Mignotte \cite{Mignotte96}, which uses linear forms in two logarithms to bound the integral coefficients of Diophantine equations of the form $ax^{m}-by^{m}=c$.
Novel advances in the general reduction of bounds arising from linear forms in 
\emph{three} logarithms have only appeared extremely recently \cite{VoutierMignotte}, demonstrating  computational finesse in \cite{BM2} and \cite{BM1}.

\subsection{Primitive divisors of Lucas--Lehmer sequences}\label{subsec:BHV}

Let $(\alpha, \beta)$ be a pair of algebraic integers. 
We say $(\alpha, \beta)$ is a Lehmer pair if $(\alpha+\beta)^2$ and $\alpha\beta$ are non-zero coprime rational integers with 
     $\alpha / \beta$  not a root of unity.
The Lehmer sequence associated to the pair $(\alpha, \beta)$ is then defined as 
\begin{equation*}
        \tilde{u}_{n}(\alpha, \beta)=
        \begin{cases}
        \frac{\alpha^{n}-\beta^{n}}{\alpha-\beta}\qquad & \text{ if $n$ is odd}\\
        \frac{\alpha^{n}-\beta^{n}}{\alpha^2-\beta^2}\qquad & \text{ if $n$ is even}.
        \end{cases}
\end{equation*}

\medskip

Similarly, $(\alpha, \beta)$ is a Lucas pair if
$\alpha+\beta$ and $ \alpha\beta$ are non-zero coprime rational integers and $\alpha / \beta$ is not a root of unity. 
The corresponding Lucas sequence is defined as
\begin{equation*}
    u_{n}(\alpha, \beta) = \frac{\alpha^n-\beta^n}{\alpha-\beta}\qquad (n=0,1,2,\dots)
\end{equation*}
which can be rewritten as 
\begin{equation*}
        u_{n}(\alpha, \beta)=
        \begin{cases}
        \tilde{u}_{n}(\alpha, \beta)\qquad & \text{ if $n$ is odd}\\
        (\alpha+\beta)\tilde{u}_{n}(\alpha, \beta)\qquad & \text{ if $n$ is even},
        \end{cases}
\end{equation*}
showing that every Lucas pair is also a Lehmer pair.

\medskip

A primitive divisor of $\tilde{u}_{n}$ is a prime divisor $p$ of $\tilde{u}_{n}$ not dividing $(\alpha^2-\beta^2)^2\tilde{u}_{1}\cdots\tilde{u}_{n-1}$ (analogously for $u_{n}$). 
The celebrated theorem of Bilu, Hanrot and Voutier \cite{BHV} is not just impressive in its own right, but lends a hand to many Diophantine applications (see \cite{GouWang}, \cite{VP21}, \cite{Pink}, \cite{GarciaPatel2020}, \cite{CIKPU}).
\begin{theorem}[Bilu, Hanrot and Voutier]\label{thm:BHV}
    Any Lucas or Lehmer sequence 
    has a primitive divisor for all $n>30$ and for all primes 
    $n>13$.
\end{theorem}

\medskip

For example, the following theorem due to Patel \cite{VP21} is proved using the aforementioned theorem of Bilu, Hanrot and Voutier in combination with a factorisation argument over imaginary quadratic fields, displaying how well exponents in exponential Diophantine equations can be bounded.
   \begin{theorem} \label{thm:BHVPatel} Let $C_{1} \geq 1$ be a squarefree integer and $C_{2}$ a positive integer. Assume that $C_{1}C_{2} \not\equiv 7 \pmod{8}$. Let $p$ be
    an odd prime for which
 \begin{equation}
    \label{eq: BHV}
     C_{1}x^{2}+C_{2}=y^{p}, \qquad \gcd(C_{1}x^{2}, C_{2}, y^{p})=1 
 \end{equation}
has a solution $(x, y) \in \bbN^{2}$ such that $xy \ne 0$. Writing $C_{1}C_{2} = cd^{2}$ with $c$ squarefree, then either
\begin{itemize}
    \item $p \leq 5$; or
    \item $p=7$ and $y=3, 5$ or $9$; or
    \item $p$ divides the class number of $\Q(\sqrt{-c})$; or
    \item $p$ divides $\left(q-\left(\frac{-c}{q}\right)\right)$, where $q$ is a prime $q \mid d$ and $q \nmid 2c$. 
\end{itemize}
\end{theorem}

\subsection*{Eliminating equations}  
Suppose we have bounded the exponent of the exponential Diophantine equation that we are interested in solving. We are now interested in concretely solving as many possibilities as we can. 
This gives us finitely many equations to resolve. We are now interested in finding the full set of solutions to a given equation or proving that it has no solutions (in other words, eliminating the equation).
In this section we summarise two of the most powerful tools used to eliminate equations: Sophie Germain's empty set criterion and Wiles' Modularity Theorem.

\subsection{Sophie Germain's empty set criterion}\label{subsec:SophieGermain}

Sophie Germain (1776--1831) made significant contributions to number theory despite being born into an era in which it was socially improper for women to practice mathematics. 
In pursuit of her dream, having adopted the identity of the former student \emph{Monsieur Auguste Antoine LeBlanc}, she enrolled at the newly opened \'{E}cole Polytechnique in Paris in 1794. 
Sophie began corresponding with Gauss after her interest in number theory piqued through a serious study of Legendre's \emph{Th\'{e}orie des nombres} and Gauss' 
\emph{Disquisitiones Arithmeticae}.
She soon became the first to provide evidence for Fermat's Last Theorem through studying residue classes of the exponent. 
Although her work was left unrecognised during her lifetime, her contributions to mathematics did receive recognition after her death. In 1882, the street \emph{Rue Sophie-Germain} was named in her honour (located in the $14^{th}$ \emph{Arrondissement} of Paris). Significantly later on, in 2003, the French Academy of Sciences officiated an annual prize to mark her legacy. 

\medskip

In one of her 14 letters to Gauss, Sophie writes (taken from \cite[pp. 358--359]{Centina}):

\medskip

\begin{displayquote}

Prenons pour exemple l’\'{e}quation m\^{e}me de Fermat qui est la plus simple
de toutes celles dont il s’agit ici.

Soit donc, $p$ \'{e}tant un nombre premier, $z^{p} = x^{p} + y^{p}$.
Je dis que si cette \'{e}quation est possible, tout nombre premier de la forme $2Np+1$ ($N$
tant un entier quelconque) pour lequel il n’y aura pas deux r\'{e}sidus) $p^{i\acute{e}me}$
puissance plac\'{e}s de suite dans la s\'{e}rie des nombres naturels divisera
n\'{e}cessairement l’un des nombres $x, y$ et $z$.

Cela est \'{e}vident, car l’\'{e}quation $z^{p} = x^{p} + y^{p}$ donne la congruance [congruence] $1 \equiv r^{sp} + r^{tp}$ dans laquelle $r$ repr\'{e}sente une racine primitive et
$s$ et $t$ des entiers.

On sait que l’\'{e}quation a une infinit\'{e} de solutions lorsque $p = 2$. Et en
effet tous les nombres, except\'{e}s 3 et 5 ont au moins deux r\'{e}sidus quarr\'{e}s
dont la diff\'{e}rence est l’unit\'{e}. Aussi dans ce cas la forme connue savoir
$h^{2} + f^{2}, 2fh, h^{2}-f^{2}$ des nombres $z[x], y$ et $z$ montre-t-elle que l’un de
ces nombres est multiple de 3 et aussi que l’un des m\^{e}mes nombres est
multiple de 5.

Il est ais\'{e} de voir que si un nombre quelconque $k$ est r\'{e}sidu puissance
$p^{i\acute{e}me}$ mod. $2Np + 1$ et qu’il y ait deux r\'{e}sidus puissance $p^{i\acute{e}me}$ m\'{e}me
mod. dont la diff\'{e}rence soit l’unit\'{e}, il y aura aussi deux r\'{e}sidus puissance $p^{i\acute{e}me}$ dont la diff\'{e}rence sera $k$.

Mais il peut arriver qu’on ait deux r\'{e}sidus $p^{i\acute{e}me}$ dont la diff\'{e}rence soit
$k$, sans que $k$ soit r\'{e}sidu $p^{i\acute{e}me}$.

Cela pos\'{e} voici l’\'{e}quation g\'{e}n\'{e}rale dont la solution me semble d\'{e}pendre
comme celle de Fermat de l’ordre des r\'{e}sidus:
$kz^{p} = x^{p} \pm y^{p}$ car d’apr\'{e}s ce que vient d’\^{e}tre dit on voit que tout nombre premier de la forme $2Np+ 1$ pour lequel deux r\'{e}sidus $p^{i\acute{e}me}$ n’ont pas le nombre $k$ pour diff\'{e}rence divise le nombre $z$ [l’un des nombres $x, y, z$]. Il suit del\'{e} que
s’il y avait un nombre infini de tels nombres l’\'{e}quation serait impossible.

\end{displayquote}

In the above correspondence, Germain essentially provides a proof for the following theorem.

\begin{theorem}[Germain]
    Let $p\geq 3$ be a prime. Let $q=2Np+1$ be an auxiliary prime, where $N$ is an integer not divisible by 3. Suppose the following hold:\\
    i) If $x^p+y^p+z^p\equiv 0$ (mod $q$) then $xyz\equiv 0$ (mod $q$);\\
    ii) $p$ is not a $p^{th}$ power residue modulo $q$.\\
    Then there are no rational solutions to the Fermat equation with exponent $p$ when $p\nmid xyz$.
\end{theorem}

\medskip

The following criterion due to Patel \cite{Patel17}, inspired by the above work of Sophie Germain, is a computational method of eliminating Fermat-type equations of signature $(p,2p,2)$. Since all computations are performed over finite fields this method is computationally highly efficient. 

\begin{theorem}
    \label{thm:emptyset}
    Let $p \ge 3$ be a prime. 
    Let $a$, $b$ and $c$ be positive integers such that 
$\gcd(a,b,c)=1$. 
Let $q=2k p+1$ be a prime that does not divide $a$. 
Define
\begin{equation*}
\mu(p,q)=\{ \eta^{2p} \; : \; \eta \in \mathds{F}_q \}
=\{0\} \cup \{ \zeta \in \mathds{F}_{q}^* \; : \; \zeta^{k}=1\}
\end{equation*}
and
\begin{equation*}
    B(p,q)=\left\{ \zeta \in \mu(p,q) \; : \; ((b \zeta+c)/a)^{2k} \in \{0,1\} \right\} \, .
\end{equation*}
If $B(p,q)=\emptyset$ then there are no solutions to $aw_{2}^{p}-bw_{1}^{2p}= cr^{2}$.
\end{theorem}

We demonstrate the computational power of this result through the following example.

\begin{example}
    In \cite{CIKPU}, the authors resolve the equation
    \begin{equation*}
(x-4r)^3 + (x-3r)^3 + \cdots + x^3 + \cdots + (x + 3r)^3+ (x + 4r)^3 = y^p
\end{equation*}
under the constraints $x,r,y,p \in \Z$, $p$ a prime, $\gcd(x,r)=1$ and $0 < r\le 10^6$. Descent, and using the divisibility conditions $2\nmid x$,  $3\nmid x$, and  $5\mid x$, lead to the ternary equation 
\begin{equation}
    \label{eq:desc5}
    3^{p-2}\cdot w_2^p-5^{2p-3}w_1^{2p}=4r^2.
\end{equation}
Once the exponent $p$ has been bounded, there are still $1, 951, 465, 447$
equations to be resolved in unknown variables $x$ and $y$.
Germain's methodology results in only $2, 821$
equations surviving after applying Theorem \ref{thm:emptyset} to \eqref{eq:desc5}.

\end{example}

\subsection{Wiles' Modularity Theorem} 
\label{subsec:modularity}

Perhaps the most fashionable (albeit computationally challenging) method of studying exponential Diophantine equations is the one pioneered by Wiles and his collaborators during the resolution of the Fermat equation (over the rational field) during the late $20$th century. We give a brief outline of the proof here.

\medskip

Suppose we have a non-trivial rational solution $(a,b,c)$ to the Fermat equation of degree $p$
\begin{equation}
    \label{Fermat}
    x^p+y^p=z^p,
\end{equation}
for prime $p\geq 5$. We can associate a semistable elliptic curve $E$ (called the \emph{Frey-Hellegouarch} curve) to this putative solution
\begin{equation*}
        E: y^2=x(x-a^p)(x+b^p),
\end{equation*}
where $E$ has conductor $N$. A key non-trivial  ingredient in this method is the modularity of $E$. The modularity of all semistable elliptic curves over $\bbQ$ was established by Wiles \cite{Wiles} and Taylor--Wiles \cite{TW95}. A few years later, the modularity of all elliptic curves over $\bbQ$ was established by Breuil-Conrad-Diamond-Taylor \cite{BCDT2001}.

\medskip

The absolute Galois group $G_{\bbQ}$ of $\bbQ$ acts on the group of $p$--torsion points on $E$.
This induces the group representation of $G_{\bbQ}$ 
\begin{equation*}
        \bar{\rho}_{E, p}:\; G_{\bbQ}\rightarrow GL_{2}(\mathds{F}_{p})
\end{equation*}
known as the mod $p$ Galois representation $\bar{\rho}_{E, p}$ associated to $E$.
By work of Mazur \cite{Mazur78}, we know that $\bar{\rho}_{E, p}$ is irreducible. 
By the modularity of $E$, there is a newform of level $N$ associated to $\bar{\rho}_{E, p}$.
The final step is to apply Ribet's level-lowering theorem \cite{Ribet90} which 
yields a newform of level 2. 
Since there are no newforms of level 2, this contradicts the existence of the putative solution.

\medskip

This initiated a movement of establishing the modularity of elliptic curves over other families of totally real number fields (see \cite{Jarvis08}, \cite{freitas2014elliptic}, \cite{cubicmod}, \cite{box2021elliptic}). 
Subsequently, the Fermat equation has been studied extensively over totally real number fields using an extension \cite[Theorem 3]{Freitas} of the modular method outlined above (see \cite{JarvisMeekin}, \cite{Freitas}, \cite{michaudjacobs2021fermats}, \cite{KJ22}).
A particular highlight is the trailblazing work in \cite{KJ22}, which is the very first instance of a full resolution to the Fermat equation over a field with degree greater than 2 (real biquadratic field specifically).

\medskip

Comparatively, less progress has been made in the direction of establishing the modularity of elliptic curves over number fields with at least one pair of complex embeddings. However, under the assumption of Serre's modularity conjecture (see \cite[Conjecture 3.1]{SengunSiksek18}), the Fermat equation has been resolved over some small imaginary quadratic fields by \c{T}urca\c{s} \cite{Turcas18, Turcas20}. 
In very recent work, Caraiani and Newton \cite{CN23} have established the modularity of all elliptic curves over imaginary quadratic fields $F$ such the Mordell-Weil group $X_{0}(15)(F)$ is finite, opening up the road for unconditional results in the direction of the resolution of Diophantine equations over imaginary quadratic fields.
\medskip

Wiles' Theorem (for which he was subsequently awarded the Abel Prize (2016)) paved the way for further advances in 
the explicit resolution of (ternary exponential) Diophantine 
equations. Some of the novelties established to date include
Kraus \cite{Kraus}, Bennett and Skinner \cite{BenSki} and Bennett, Vatsal and Yazdani \cite{BenVatYaz} which have generalised the construction of an elliptic curve out of hypothetical solutions of generalised Fermat equations of signatures $(p, p, p)$, $(p, p, 2)$ and $(p, p, 3)$ respectively. For these three cases, the association process is very detailed. Moreover, thanks to the well-written and explicit survey \cite{Siksek}, using Wiles' modular method is more accessible and the enthusiastic Diophantine solver may employ these tactics to arrive at a contradiction.

\subsection*{Determining all integral solutions}
Ideally, we would like to determine all of the integral points of a certain Diophantine equation from the offset. However, current tools to achieve this are either not applicable in general cases, or otherwise computationally expensive and infeasible. Once we manage to sieve out impossibilities, we should hopefully be left with a handful of equations that in the end require a resolution. In this section, we distinguish 
between the common types of equations encountered, once we have finished sieving.

\subsection{Integral Points on genus 1 curves}

\label{subsec:IntPoints}
In 1929, Siegel \cite{Siegel} announced the following groundbreaking result.

\begin{theorem}[Siegel]\label{thm:Siegel}
Let $C$ be a smooth affine algebraic curve of genus $g$ defined over a number field $K$.
Then there are only finitely many $K$-integral points on $C$, provided that $g$ is strictly positive.
\end{theorem}

Siegel's Theorem is the first major result within Diophantine equations that depends only on the genus of the curve and not directly on the form of the equations. Setting $g=1$ and $K=\Q$ in Siegel's Theorem tells us that there are finitely many integer points on a given elliptic curve. Siegel used Thue's methods from Diophantine approximation to prove this result, however, since Thue's methods are ineffective, so is the proof of Siegel's Theorem. Not all is lost, as Theorem~\ref{thm:Baker} (due to Baker) 
is applicable in this case, which yields effective results that are absolutely key to perform explicit computations. We demonstrate this with some examples.

\medskip

Let us consider the elliptic curve $E/\bbQ$ with Cremona label \texttt{11a1} with minimal Weierstrass equation $E: y^2+y=x^3-x^2-10x-20$.
The Mordell-Weil group of $E/\bbQ$ is 
$E(\bbQ)=\bbZ/5\bbZ=\{0_{E}, (5,5), (5,-6), (16,60), (16,-61)\}$.
Since $E$ has finitely many rational points, we automatically know that $E$ has finitely many integral points.
Let us now consider the elliptic curve $E^{\prime}/\bbQ$ with Cremona label \texttt{37a1} with minimal Weierstrass equation $E^{\prime}: y^2+y=x^3-x$.
The Mordell-Weil group of $E^{\prime}/\bbQ$ is given by
$E^{\prime}(\bbQ)=\bbZ=\langle (0,0)\rangle$.
In other words, $E^{\prime}$ has infinitely many rational points. 
From this information, we cannot deduce the state of the integral points on $E^{\prime}$ straight away.
The determination of the Mordell-Weil group of an elliptic curve over $\bbQ$ is generally a difficult problem, but once this has been completed, then the set of integral points can be computed with ease. 
Suppose an elliptic curve $E/\bbQ$ has rank $r$. 
We can write any rational (and indeed integral) point as $P=Q+a_{1}R_{1}+\cdots+a_{r}R_{r}$ where $a_{i}\in\bbZ$, $Q$ is torsion point on $E$ and $R_{1},\cdots,R_{r}$ generate $E(\bbQ)$ modulo torsion. If $P$ is an integral point then by Theorem~\ref{thm:Siegel} there are only finitely many possibilities for $P$.
Effective bounds on $|a_{i}|$ is due to Baker, which can be reduced 
significantly further by lattice reduction algorithms (see \cite[Section 2.3.5, Section 8.7]{cohenVI} for details). 
Once a reasonable bound on the $|a_{i}|$ is obtained, it becomes feasible to completely determine the set of integral points on $E$.
We note that this method has been well-implemented computationally.
For example, computing the integral points on the elliptic curve $E^{\prime}$ on the online \texttt{MAGMA} calculator (\href{http://magma.maths.usyd.edu.au/calc/}{http://magma.maths.usyd.edu.au/calc/}) via the \texttt{IntegralPoints} command \cite{Magma} took 0.520 seconds. 

\subsection{Integral points on higher genus curves}\label{subsec:Chabauty}
 
 In many Diophantine problems it is common to run into hyperelliptic curves, for example, after performing a descent step. It is possible to solve some of these equations by using the method of Chabauty via a \texttt{MAGMA} implementation \cite{Magma}. The aim of Chabauty's method is to recover the set of rational points $C(\Q)$ of a smooth projective curve $C / \Q$ from its Jacobian $J(C)$ through the Abel-Jacobi map when the rank of $J(C)$ is not ``too big''. We refer the reader to \cite{Siksek2015Chabauty} or \cite{MullerChabauty} for a comprehensive survey of the methodology and we proceed to summarise the main ideas below. 

\medskip

Let $C/\Q$ be a curve of genus $g$ and let $J(C)$ denote its Jacobian. Assume that we have a known rational point $P_{0} \in C(\Q)$; the Abel-Jacobi embedding can be written as 
$$\iota: C \rightarrow J(C), \qquad P \mapsto [P-P_{0}].$$
Consider $\Omega_{C} := H^{0}(C, \Omega_{1})$ the $g$-dimensional vector space of regular differentials of the curve $C$. Given a prime $p \in \bbN$, we have the pairing
$$ \Omega_{C} \times J(\Q_{p}) \rightarrow \Q_{p} $$
defined by 
$$ \langle \omega, [\sum_{i} P_{i}-Q_{i}]\rangle := \sum \int_{Q_{i}}^{P_{i}} \omega. $$
If the divisor  $[\sum_{i} P_{i}-Q_{i}]$ is a torsion element of $J(\Q_{p})$, then we have $\sum \int_{Q_{i}}^{P_{i}} \omega = 0$. Thus $J(\Q)_{\text{tor}} \subset J(\Q_{p})_{\text{tor}} $ lies in the right kernel of the pairing $\langle, \rangle$. Thus, if we have $D_{1}, \ldots, D_{r}$ a basis for the free part of $J(\Q_{p})$, then the space $J(\Q)^{\perp}$ of differentials $\omega$ annihilated by all elements $D$ in $J(\Q)$ is determined by the ones annihilated by $D_{1}, \ldots, D_{r}$. That is, $J(\Q)^{\perp}$ is given by the equations 
$$ \int_{D_{i}} \omega = 0, \qquad i=1, \ldots, r.$$
It follows then that $\dim J(\Q)^{\perp} > g-r$. If the condition $r<g$ (known as the Chabauty condition) holds, then $\dim J(\Q)^{\perp} \geq 1$ and there exists a differential $\omega_{0} \in \dim J(\Q)^{\perp}$. Then, the points $C(\Q)$ are determined by the equation $$ \int_{\iota{P}} \omega_{0} = 0.$$
This equation with a finite number of solutions can be solved by using $p$-adic power series. Finally, the set of rational points $C(\Q)$ can be determined using a Mordell-Weil sieve. The combination of these techniques is currently implemented in \texttt{MAGMA}  \cite{Magma} for hyperelliptic curves. 
The implementation is often capable of determining $C(\Q)$.

\begin{example} In \cite{CIKPU}, 
the authors resolve the equation
 \begin{equation*}
(x-4r)^3 + (x-3r)^3 + \cdots + x^3 + \cdots + (x + 3r)^3+ (x + 4r)^3 = y^p
\end{equation*}
under the constraints $x,r,y,p \in \Z$, $p$ a prime, $\gcd(x,r)=1$ and $0 < r\le 10^6$. 
For $p=5$, descent and the divisibility conditions $3 \nmid x$ and $4 \mid x$ and $5 \nmid x$ lead to the ternary equation
$3^{3}w_{2}^{5}-2^{4}w_{1}^{10}=5r^{2}$.
Letting $X=w_{2}/w_{1}^{2}$ and  $Y=5r$ yields the genus 2 curve
\[
C:\; Y^2=3^{3}\cdot 5\cdot X^{5}-2^{4}\cdot 5.
\]
The Jacobian of $C$ has rank 1 over $\Q$. Thus, the Chabauty condition is satisfied and the Chabauty implementation in \texttt{MAGMA} \cite{Magma} is used to find that
\[
C(\Q)=\{\infty\}.
\]

\end{example}

\subsection{Irreducible bivariate forms}\label{subsec:Thue}

A Thue equation of degree $n\geq 3$ is a Diophantine equation of the form
\begin{equation}
    \label{eq:Thue}
    f(x,y)=C
\end{equation}
where $f$ is an irreducible degree $n$ polynomial with integer coefficients  and $C$ is a known integer.
It is well-known, due to a non-effective proof of Thue \cite{Thue09}, that there are finitely many integral solutions to \eqref{eq:Thue}. 

\medskip

Effectivity was first achieved via Baker's theory of linear forms in logarithms.
There are now readily available algorithms (implemented in \texttt{MAGMA} \cite{Magma}) that solve \eqref{eq:Thue} for small values of $n$. 
For a detailed survey, we refer the reader to \cite{Waldschmidt20}. 

\begin{example}
    In \cite{CIKPU}, the authors resolve the equation
    \begin{equation*}
    (x-4r)^3 + (x-3r)^3 + \cdots + x^3 + \cdots + (x + 3r)^3+ (x + 4r)^3 = y^p
    \end{equation*}
under the constraints $x,r,y,p \in \Z$, $p$ a prime, $\gcd(x,r)=1$ and $0 < r\le 10^6$.
Within each descent case, there are a small number of equations that remain once we have sieved out possibilities. These are usually dealt with using a Thue equation solver. Following previous approaches taken in \cite{GarciaPatel2019}, \cite{Garcia2019} and \cite{GarciaPatel2020} proved to be 
futile as 
the authors encountered Thue equations with 
 very large exponents and coefficients. These were
 (miraculously) eliminated via an application of the methods presented in Section~\ref{subsec:BHV}. Here is an example of a Thue equation that was bypassed: 
\[
X^{19} - 48536986868323784437961876392364501953125(Y^2)^{19} = 274877906944.
\]
We note that this particular Thue equation can (just about) be solved using \texttt{MAGMA}'s \cite{Magma}
inbuilt Thue solver, under GRH assumption. 
\end{example}

\section{Perfect powers as sums of like powers}

In this section, we  give an overview of the results to date on various subfamilies of \eqref{eq:table1};
\begin{equation*} 
x^{k}+(x+r)^{k}+ \cdots + (x+(d-1)r)^{k} =y^{n} 
\qquad d, k, n, x, y, \in \Z,  \quad r \geq 1; \quad d, k, n\geq 2 
\end{equation*}
with parameters $(k, d, r, n)$. We present tabular summaries 
at the end of each subsection.

\subsection{Sums of squares}

Let's consider the problem of determining perfect powers that are sums of $d$ consecutive squares, i.e., we consider the following subfamily of  \eqref{eq:table1}: 
\begin{equation}\label{eq:conssquares} 
x^{2}+(x+1)^{2}+ \cdots + (x+d-1)^{2} =y^{n} \qquad d, n, x, y, \in \Z,
\quad d,n\geq 2. 
\end{equation}

\medskip

The case $d=2$ was first studied by Cohn \cite{cohn_1996} in 1996 using elementary number theory techniques in conjunction with the Fibonacci sequence.  
It is well-known that \eqref{eq:conssquares} with exponent $n=2$ possesses
infinitely many solutions. For example, \eqref{eq:conssquares} in the case $n=2$ and $d=2$ admits the following infinite family of solutions (\cite{patel2017powers}),
 \[
     2x+2+y\sqrt{2}=\pm(1+\sqrt{2})^{2a+1},\; a\in\bbZ.
 \]

\medskip

In \cite{BreStrTza}, the authors study \eqref{eq:conssquares} also with $n=2$, however, they fix values of $x$ and let $d$ vary. In this manner, they are able to re-frame the problem to one of determining integral points on certain elliptic (genus 1) curves. Letting  $1\leq x\leq 100$, Bremner, Stroeker and Tzanakis found that the principal difficulty arose in determining the torsion-free part of the curve. 
The authors subsequently used the method of linear forms in logarithms to obtain a sensible bound on the heights of integral points on the curve,
allowing them to feasibly search for all such points.

\medskip

After the passing of a couple of decades, Zhang and Bai \cite{ZhangBai2013} provided a resolution of \eqref{eq:conssquares} when \ $d=x$ using work of Bazs\'{o}, B\'{e}rczes, Gy\H{o}ry, and Pint\'{e}r \cite{BazBerGyoPint} on the equation $Ax^{n} - By^{n} = C$. A year later,
Zhang \cite{zhang2014diophantine}, using algebraic number theory, provides a complete resolution to the problem of determining when sums of three consecutive squares forms a perfect power; he proves that this cannot happen.

\medskip

In an extension to this problem, one may consider \eqref{eq:conssquares}, where the terms are taken in (bounded) arithmetic progression. Application of work of Bilu, Hanrot and Voutier (see Section~\ref{subsec:BHV}) leads to the resolution of further general subfamilies, 
see \cite{van18sq, koutsianas2018perfect, kundu2021perfect, LeSoy2022}.

\medskip

The first glimmer of density results can be found in 
\cite{ZhangBai2013}. Theorem 1.2 in \cite{ZhangBai2013}
states that there are infinitely many values of $d$ such that
the sum of $d$ consecutive squares will never be a perfect power. A straightforward application of Dirichlet's Theorem gives the desired density result. We further explore density approaches in Section~\ref{subsec:higherpow}.

\begin{table}[htbp!]
    \centering
    \begin{tabular}{ccc}
    \specialrule{.2em}{.1em}{.1em}  
         $(k, d, r, n)$ & Reference(s) & Method(s) used \\
         \specialrule{.2em}{.1em}{.1em} 
         $(2, 2, 1, n \geq 3)$ & Cohn \cite{cohn_1996} & elementary  \\ \hline
         \makecell{$(2, d \geq 2, 1, 2)$ \\ for $1 \leq x \leq 100$} & \makecell{Bremner, Stroeker \\ and Tzanakis \cite{BreStrTza}} & \ref{subsec:linlogs}, \ref{subsec:IntPoints} \\ \hline
         $(2, d=x, 1, \text{ odd } n > 2)$ & Zhang and Bai \cite{ZhangBai2013}  & applied \cite{BazBerGyoPint},  \ref{subsec:IntPoints}  \\
         \hline
         $(2, 3, 1, n \geq 2)$ & Zhang \cite{zhang2014diophantine} & algebraic number theory\\   \hline
                 $(2, 2 \leq d \leq 10, 1, n \geq 3)$ & Patel \cite{van18sq} & elementary, \ref{subsec:BHV} \\ \hline
         $(2, 3, 1 \leq r \leq 10^6, n \geq 2)$ & Koutsianas and Patel \cite{koutsianas2018perfect} & \ref{subsec:BHV} \\  \hline
         $(2, 2\leq d \leq 10, 1 \leq r \leq 10^6, n \geq 2)$ & Kundu and Patel \cite{kundu2021perfect} & \ref{subsec:BHV} \\ 
         \specialrule{.2em}{.1em}{.1em} 
    \end{tabular}
    \caption{Progress towards equation \eqref{eq:table1} with $k=2$.}
\end{table}

\subsection*{Open Problems}
\begin{itemize}
\item Solve \eqref{eq:conssquares} for all values of $d$ and $n\geq 3$.
\item Solve \eqref{eq:conssquares} with fixed value of $d$, and unrestricted arithmetic progression for $n \geq 3$.
\item As of now, the problem of perfect powers that are sums of squares 
has not been particularly amenable to the modular approach (as well as generally being ungovernable by the approaches outlined in \ref{subsec:linlogs}, \ref{subsec:SophieGermain}). Could a new take on the modular approach hope for some resolution?
\end{itemize}

\medskip

\subsection{Sums of cubes}

Whilst one can easily make some progress on finding 
perfect powers that are sums of squares via  
some elementary algebraic number theory, the question of 
finding perfect powers as sums of cubes 
enters a completely different realm of 
number theory.
\begin{equation} 
    \label{eq:conscubes}
        x^{3}+(x+1)^{3}+ \cdots (x+d-1)^{3} = y^{n}, 
        \qquad x, y, d , n \in \Z, 
        \quad d,n \geq 2.
\end{equation}

\medskip

Finding squares and cubes that are a sums of a fixed number of consecutive cubes (\cite{Cassels}, \cite{Uchiyama}, \cite{Stroeker1995}) nowadays boils down to finding integral points on certain genus 1 curves \cite{bennett2017perfect}.

\medskip

Finding perfect powers (beyond squares and cubes) that are sums of consecutive cubes 
has only been touched on recently, due to contemporary techniques 
required in order to reach a full resolution. Zhang \cite{zhang2014diophantine}, building upon Cassels' work \cite{Cassels} by a clever application of the innovative and impressive work on \emph{the generalised Fermat equation} presented in \cite{BGMP2006} and \cite{BenVatYaz},  determines \emph{all} perfect powers that arise as a sum of three consecutive cubes. Zhang recovers Euler's identity, in addition to the solutions provided by Cassels, 
and concludes that there are no further solutions for all $n >1$. 

\medskip

The first approach unifying many theories in the 
world of explicit number theory was given by 
Bennett, Patel and Siksek \cite{bennett2017perfect}, who provide a 
complete list of solutions for perfect powers that are sums of at most $50$ consecutive cubes. The methods employed include 
those summarised in Sections~\ref{subsec:linlogs}, \ref{subsec:IntPoints}, \ref{subsec:SophieGermain}, \ref{subsec:modularity} and  \ref{subsec:Thue}.
 In theory it would be possible, via existing techniques, to reach satisfactory conclusions for \emph{any} value of $d$, however,
current computational techniques used would need to be enhanced in order for the computations to be completed within a reasonable time.

\medskip

Extending \cite{bennett2017perfect} to the study of perfect powers
that arise from cubes of an arithmetic progression, we note the following papers: \cite{GarciaPatel2019}, \cite{Garcia2019}, \cite{GarciaPatel2020} and \cite{CIKPU}. The exact number of cubes required is specified, and the length of the arithmetic progression is bounded (usually by $10^6$). As the number of cubes increases, complexities arise 
due to the synthesis of vastly differing theories. Whilst usually this is the biggest hindrance, the authors in \cite{CIKPU} have managed to synergise various theories to drastically increase efficiency and find all perfect powers that are sums of nine cubes with bounded arithmetic progression. Disappointingly, for $p\geq 5$, no solutions were found! In order to reach contradictions for such a large volume of equations,
refining existing methodologies (\ref{subsec:linlogs}, \ref{subsec:SophieGermain}, \ref{subsec:BHV}) was crucial.

\medskip

We reveal the state of the art here (noting that it is natural to take $n=p$ to be a prime, as we can absorb excess powers into $y$, see  \cite{CIKPU} for full details):
\begin{theorem}[Coppola, Curc\'o-Iranzo, Khawaja, Patel and \"{U}lkem]
    For a prime $p\geq 5$, there are no non-trivial integer solutions to 
    \begin{equation}
        (x-4r)^3+(x-3r)^3+(x-2r)^3+(x-r)^3+x^3+(x+r)^3+(x+2r)^3+(x+3r)^3+(x+4r)^3=y^p
    \end{equation}
    for coprime $x,r$ and $1 \leq r\leq 10^{6}$, where the trivial solutions satisfy $xy=0$.
\end{theorem}
Ideally, one would like to resolve \eqref{eq:conscubes} for any arithmetic progression. We note that current modularity methods are unable to achieve this and completely new ideas are required.

\medskip

As of now, there exists one complete resolution of equation~\eqref{eq:conscubes} with  unrestricted arithmetic progression. Edis \cite{EdisThesis} solves the problem of perfect powers that are sums of three cubes in \emph{any} arithmetic progression, 
proving that no prime power greater than 7 equals the sum of three cubes in an arithmetic progression, thus engulfing the results presented in \cite{GarciaPatel2019}.
Edis uses modularity techniques and by twisting the Frey-Hellegouarch curves, he is able to elegantly reach a contradiction. 

\begin{table}[htbp!]
    \centering
    \begin{tabular}{ccc}
    \specialrule{.2em}{.1em}{.1em}  
         $(k, d, r, n)$ & Reference(s) & Method(s) used \\
         \specialrule{.2em}{.1em}{.1em} 
         $(3, 3, 1, 2)$ & Cassels \cite{Cassels}, Uchiyama \cite{Uchiyama} & elementary  \\ \hline
                  $(3, 3, 1, n \geq 2)$ & Zhang \cite{zhang2014diophantine} &  applied \cite{Cassels}, \cite{BGMP2006} and \cite{BenVatYaz} \\  \hline
                  $(3, 2 \leq d \leq 50, 1, 2)$ & Stroeker \cite{Stroeker1995} & \ref{subsec:linlogs} \\ \hline
         $(3, 2 \leq d \leq 50, 1, n\geq 2)$ & Bennett, Patel and Siksek \cite{bennett2017perfect}  &
         \ref{subsec:linlogs},
         \ref{subsec:IntPoints},
         \ref{subsec:SophieGermain},
         \ref{subsec:modularity},
         \ref{subsec:Thue}\\ 
         \hline 
         $(3, 3, 1 \leq r \leq 10^6, \text{prime } p\ge 2)$ & Arg\'aez-Garc\'ia and Patel \cite{GarciaPatel2019} &  
         \ref{subsec:linlogs},
         \ref{subsec:IntPoints},
         \ref{subsec:SophieGermain},
         \ref{subsec:modularity},
         \ref{subsec:Chabauty},
         \ref{subsec:Thue} \\ \hline    
         $(3, 3, r \ge 1, \text{prime } p\ge 7)$ & Edis \cite{EdisThesis} &  
         \ref{subsec:modularity}\\ \hline
         $(3, 5, 1 \leq r \leq  10^6, \text{prime } p\ge 5)$ &  Arg\'aez-Garc\'ia \cite{Garcia2019} & 
         \ref{subsec:linlogs},
         \ref{subsec:BHV},
         \ref{subsec:SophieGermain},
         \ref{subsec:Chabauty},
         \ref{subsec:Thue}\\ \hline
         $(3, 7, 1 \leq r \leq  10^6, \text{prime } p\ge 5)$ & Arg\'aez-Garc\'ia and Patel \cite{GarciaPatel2020} & 
         \ref{subsec:linlogs},
         \ref{subsec:BHV},
         \ref{subsec:SophieGermain},
         \ref{subsec:Chabauty},
         \ref{subsec:Thue}\\ \hline
         $(3, 9, 1\leq r \leq 10^6, \text{prime } p\ge 2)$ & \makecell{ Coppola,  Curc\'o-Iranzo, Khawaja, \\
          Patel and \"{U}lkem \cite{CIKPU}}
        & 
         \ref{subsec:linlogs},
         \ref{subsec:IntPoints},
         \ref{subsec:BHV},
         \ref{subsec:SophieGermain},
         \ref{subsec:Chabauty},
         \ref{subsec:Thue}\\ 
          \specialrule{.2em}{.1em}{.1em} 
             \end{tabular}
    \caption{Progress towards equation \eqref{eq:table1} with $k=3$.}
\end{table}

\subsection*{Open Problems}
\begin{itemize}
\item Study families of (or sporadic) solutions for \eqref{eq:conscubes} for any $d$ and any arithmetic progression with $n=3$.
\item Solve \eqref{eq:conscubes} for all values of $d$ and $n\geq 5$.
\item Solve \eqref{eq:conscubes} with fixed value of $d$, and unrestricted arithmetic progression for $n \geq 5$.
\item Density results? Or an average statement?
\end{itemize}

\subsection{Beyond sums of squares or cubes}
\label{subsec:higherpow}
A natural continuation of the preceding subsections is to consider \eqref{eq:table1} with exponents $k$ greater than $3$.
As one can imagine, the literature now begins to rapidly thin out!
We begin with the work of Zhang once more. In \cite{zhang2014diophantine},
it is shown that sums of three consecutive fourth powers is never a perfect power. 
In this instance, the author employs the modular method (by hand!) following the recipe provided in \cite{BenSki} which attaches Frey-Hellegouarch curves to generalised Fermat-type equations of signature $(p,p,2)$, along with considering integral points on certain elliptic curves, as well as solving a handful of Thue equations to reach a full resolution (Section~\ref{subsec:modularity}, \ref{subsec:IntPoints}, \ref{subsec:Thue}).

\medskip

Bennett, Patel and Siksek \cite{Bennett2016} further this work 
and show that sums of three consecutive fifth or sums of three consecutive sixth powers is never a perfect power. In order to reach such a resolution, the modular method is significantly 
improved upon to handle larger dimensional newform spaces, notable via a clever consideration of their Hecke polynomials.
The authors employ a \emph{Multi-Frey-Hellegouarch} approach, pioneered by Bugeaud, Luca, Mignotte and Siksek \cite{BugMigSik}. The crux of the idea is to associate two (or more) different curves to the Diophantine equation, which results in performing computations within two different newform spaces that must still be compatible via modularity of the associated Galois representations. This increases the chances of reaching contradictions. We emphasise that the results in \cite{Bennett2016} for exponent $k=6$ is only possible as the authors attach Frey-Hellegouarch curves to certain binary cubic forms, as pioneered in \cite{BenDahmen}.

\medskip

The echoing question of higher powers in arithmetic progression is once again considered first by Zhang  \cite{zhang2017diophantine4} 
with partial results when considering fourth powers that depend upon the arithmetic progression.
Remarkably, van Langen \cite{van2021sum} provides a 
full resolution to the question of perfect powers that are sums of three fourth powers in \emph{any} arithmetic progression. A modular
approach is undertaken, this time with a careful analysis of certain $\Q$-curves.

\medskip

Building upon \cite{Bennett2016}, Bennett and Koutsianas \cite{BennettKoutsianas} completely solve the question of perfect powers expressible as sums of three fifth powers in arithmetic progression (superseding earlier work \cite{DDKT}). 
They dealt with small exponents using the method of Chabauty over number fields.

\medskip

To date, only a single result exists towards studying the \emph{entire} Diophantine family \eqref{eq:table1}. 
Patel and Siksek \cite{patel2017powers} prove that in the case of
even exponents $k$, and any arithmetic progression $r$, 
solutions are rare. 

\begin{theorem}[Patel and Siksek]
    \label{thm: PS17}
    Let $k\geq 2$ be an even integer, $r\in\bbZ_{> 0 }$.
    Let $\mathcal{A}_{k,r}$ denote the set of integers $d\geq 2$ such that there is a solution to \eqref{eq:table1} for $n\geq 2$. Then the set $\mathcal{A}_{k,r}$ has natural density 0.
\end{theorem}

The proof relies on the study of the Galois group and the Newton polygon of Bernoulli polynomials. Bernoulli polynomials are conjectured to be an irreducible family of polynomials, however, 
at this moment, we seem to be quite far away from a resolution to 
this long-standing conjecture.
We further remark that  we can not hope for such a strong density result for odd values of $k$. This is due to the existence of the solution $(x,y,n)=(r(1-d)/2,0,n)$, where $d$ is any odd integer.

\begin{table}[htbp!]
    \centering
    \begin{tabular}{ccc}
    \specialrule{.2em}{.1em}{.1em}  
         $(k, d, r, n)$ & Reference(s) & Method(s) used \\
         \specialrule{.2em}{.1em}{.1em} 
          $(4, 3, 1, n \geq 2)$ & Zhang \cite{zhang2014diophantine} &  \ref{subsec:IntPoints}, \ref{subsec:modularity}, \ref{subsec:Thue} \\  \hline
          $(4, 3, 1 \leq r \leq 6, n \geq 2)$& Zhang \cite{zhang2017diophantine4} & \ref{subsec:IntPoints}, \ref{subsec:modularity}, \ref{subsec:Thue}\\ \hline
           $(4, 3, r \ge 1, n\geq 2)$ & Van Langen  \cite{van2021sum} & \ref{subsec:modularity}, analysis of $\Q$-curves\\ \hline
                    $(5 \leq k \leq 6, 3, 1, n \geq 2)$ & Bennett, Patel and Siksek \cite{Bennett2016}  & \makecell{\ref{subsec:modularity}, \ref{subsec:Thue}, \ref{subsec:IntPoints}, 
         \\ Hecke polynomials} \\ \hline
         $(5, 3, r \ge 1, n\geq 2)$ & Bennett and Koutsianas \cite{BennettKoutsianas} & \ref{subsec:modularity}, \ref{subsec:Chabauty} \\ \hline
         \makecell{``For even $k$, and $n \ge 2$, \\ equation \eqref{eq:table1} has  no solutions \\ for almost all $d\geq 2$''} & Patel and Siksek \cite{patel2017powers} & 
         \makecell{Bernoulli polynomial,  and \\ its Galois Group and \\  Newton Polygon} \\ 
         \specialrule{.2em}{.1em}{.1em} 
    \end{tabular}
    \caption{Progress towards equation \eqref{eq:table1} for $k\geq 4$.}
\end{table}

 \subsection*{The asymptotic approach: a way forward?}
 In 2014, Freitas and Siksek initiated the study of viewing the Fermat equation asymptotically. As in \cite{freitas_siksek_2015}, we say \emph{Asymptotic Fermat's Last Theorem} holds over a (totally real) number field  $K$ if there is a constant $B_{K}$ such that for all $p>B_{K}$ there are no non-trivial solutions to the Fermat equation over the ring of integers of $K$, namely, $\mathcal{O}_K$. 

\medskip

The asymptotic approach begins, in the same way as the modular approach, by attaching a Frey curve $E$ to a putative non-trivial solution to the Fermat equation. 
If $E$ is modular over $K$, then $E$ corresponds to some modular form. 
Level-lowering produces a Hilbert newform $\mathfrak{f}$ of level $\mathcal{N}_{p}$ such that $\bar{\rho}_{E, p}\sim \bar{\rho}_{\mathfrak{f}, \bar{\omega}}$. 
In a diversion to the modular approach, Freitas and Siksek then use a Eichler-Shimura type assumption to construct an elliptic curve $E_{1}: y^2=(x-e_{1})(x-e_{2})(x-e_{3})$ of conductor $\mathcal{N}_{p}$, corresponding to $\mathfrak{f}$.
Let $S$ denote the set of primes of bad reduction for $E_{1}$.
Writing $\lambda=(e_{1}-e_{2})/(e_{1}-e_{3}), \mu=(e_{2}-e_{3})(e_{1}-e_{3})$ gives a solution to the $S$--unit equation $\lambda+\mu=1$, where $S$ is the set of primes dividing $2\mathcal{N}_{p}$.
The idea is then to study solutions to the Fermat equation over $K$ by studying solutions to the  $S$--unit equation over $K$. This method has since been extended to study other Fermat-type equations (see \cite{Deconinck15}, \cite{SengunSiksek18}, \cite{Kraus19}, \cite{KaraOzman20}, \cite{FREITAS2020106964},  \cite{IsikKaraOzman20}, \cite{isik_kara_ozman_2022}, \cite{mocanu222} and \cite{mocanu221}).

\subsection*{Open Problems}
\begin{itemize}
\item Can we find all solutions to \eqref{eq:table1} for some fixed $d$ and $r$ with $n>3$?
\item Solve \eqref{eq:table1} for higher exponents $k \geq 4$ and larger length $d$ of consecutive terms. 
\item Solve \eqref{eq:table1} for higher exponents $k\geq 4$ with fixed  $d$ and unrestricted arithmetic progression. 
\item Consider perfect powers that are sums of differing powers as studied in  \cite{zhang2017diophantine}.
\item Following \cite{zhang2017diophantine}, is there any hope to solve this for differing powers taken in an arithmetic progression? What about \emph{all} arithmetic progressions?
\item Density results for odd exponents, modulo trivial solutions? Or an average statement?
\item Proving irreducibility of the entire family of Bernoulli polynomials with even degree.
\item Proving irreducibility of the entire family of Bernoulli polynomials modulo the factors arising from the well-known rational roots for all odd degree.
\end{itemize}

\section{Perfect powers as sums of consecutive integers starting from 1}
 In this section, we present an overview of the results to date on various subfamilies of \eqref{eq:table2}:
\begin{equation*} 
s(1^k+2^k+ \cdots + x^k) + r(x) =y^{n} 
\qquad k, n, x, y, s \in \Z, , r(x) \in \Z[x],
\quad k,n \geq 2
\end{equation*}
with parameters $(k, n, s, r(x))$. Notice that when $s=1$ and $r(x)=0$ this equation gives a subfamily of equation \eqref{eq:table1}. A note of great importance is that the methodology and results obtained for this subfamily differ vastly 
 from those given in the previous section. We summarise key contributions at the end of each subsection in Tables~\ref{tab:EM} and \ref{tab:Table2}, as well as listing open problems.

\medskip

 This story  begins with Lucas' Problem 1180 \cite{Lucas1180}:
\begin{displayquote} 
    Une pile de boulets \`a base car\'ee ne contient un nombre de boulets \'egal au carr\'e d’un nombre entier que lorsqu’elle en contient  vingt-quatre sur le c\^ot\'e de la base.
\end{displayquote}

\medskip

In this assertion, E. Lucas claimed that the only (positive integer) solutions to the equation 
\begin{equation*}
    1^2+\cdots+x^2=y^2 
\end{equation*}
are $(x, y) \in \{(1, 1), (24, 70)\}$. Attempts to prove this claim by himself \cite{LucasSol1180} and Moret-Blanc \cite{MoretBlanc} turned out to be both defective, as they contained errors in the argumentation. It was only in 1918 that a correct proof materialised. 

\medskip

Watson \cite{Watson} was the first to do so, but many others have given alternative proofs of the theorem; see \cite{Ljunggren}, \cite{Ma1985}, \cite{MaDiophantine}, \cite{Cucurezeanu}, \cite{Anglin} and \cite{BennettPyramid}.

\subsection{Erd{\H o}s-Moser Conjecture}
In the early 1950's, Erd\H{o}s and Moser studied the more general equation
\begin{equation}
    \label{eq:erdmos}
S_k(m):= 1^k + 2^k + \cdots + (m-1)^k = m^k.
\end{equation} 
It is easy to see that for $k=1$ we can rewrite the left hand side of \eqref{eq:erdmos} as $m(m-1)/2$ and find the unique trivial solution $m=2$. However, as Erd\H{o}s noticed in a letter to Moser, taking $k\geq 2$, equation \eqref{eq:erdmos} seems to have not a single solution. 

\medskip

Most Diophantine problems (of interest) seem to be easy to state, yet notoriously 
difficult to solve; naturally this is the case for equation~\eqref{eq:erdmos}, as it still remains conjectural.

\medskip

\begin{conjecture}[Erd\H{o}s-Moser]\label{conj:ErdosMoser} 
The equation \eqref{eq:erdmos}
only has one integral solution when $k=1$ and $m=2$ and no others.
\end{conjecture}

\medskip

There exist (partial) results that should convince the reader to stay positive about Conjecture \ref{conj:ErdosMoser} being true.  For odd $k$, it is easy to see that the conjecture is indeed true. By induction, one can see that for such $k$ the sum $S_{k}(m)$ needs to be a multiple of the integer $m(m-1)/2$; provided that $m>2$, if there would be a solution to \eqref{conj:ErdosMoser} this would contradict the coprimality of $m$ and $m-1$. Hence $k$ has to be even.
Moser himself \cite{Moser1953}, assuming that a solution exists to \eqref{eq:erdmos}, with $k\ge 2$ even, ingeniously provided a proof for a large lower bound on $m$.

\begin{theorem}[Moser, 1953]\label{thm:Moser}
If $(m,k)$ is a solution to \eqref{eq:erdmos}, with $k\ge 2$, then $m>10^{{10}^{6}}$.
\end{theorem}

Moser's proof unfolds from four simple identities (or mathemagical rabbits) that putative solutions of \eqref{eq:erdmos} need to satisfy. Moree \cite{MoreeRabbits} gave a top hat to these magical identities by noticing that they can be pulled up from a theorem on harmonic series \textit{\`{a} la} von Staudt-Clausen's Theorem. 
By doing some juggling with the four mathemagical rabbits of Moser, one can arrive to the conclusion that any solution $m$ would need to satisfy 
$$\sum_{p\mid M} \frac{1}{p} + \frac{1}{m-1} + \frac{2}{m+1} + \frac{2}{2m-1} + \frac{4}{2m+1} > 3.16 $$
where $M$ is a squarefree integer that can be constructed from $m$. 
Now the lower bound on $m$ follows from looking for the smallest such $M$ such that the value of prime harmonic series -- which is a divergent series -- is big enough. Moser estimated this value, thus obtaining an ``estimated value" for the lower bound on $m$. 

\medskip

With the arrival of the new millennia and the emergence of modern techniques, Butske, Jaje and Mayernik \cite{BuJaMay} managed to refine Moser's lower bound by explicitly computing the value of $m$ derived from the original outlined proof. The improvements situated $m$ just above $1.485\cdot10^{9321155}$.

\begin{theorem}[Butske, Jaje and Mayernik, 2000]\label{thm:BJM}
If $(m,k)$ is a solution to \eqref{eq:erdmos}, with $k\ge 2$, then $m>10^{9\cdot {10}^{6}}$.
\end{theorem}

Notice that Butske, Jaje and Mayernik \cite{BuJaMay} are just shy of the $10^{{10}^{7}}$ mark! The authors remark 
``... the authors hope that new insights will eventually make it possible to reach the more natural benchmark $10^{{10}^{7}}$..."

\medskip

Readers, fear not, as the change of millennia also brought with it completely novel methods: Gallot, Moree and Zudilin used continued fractions digits calculation of $\log 2/N$ for an  integer $N$ and managed to obtain significant improvements on the lower bound $m$. 

\begin{theorem}[Gallot, Moree and Zudilin, 2011]\label{thm:GMZ}
If $(m,k)$ is a solution to \eqref{eq:erdmos}, with $k\ge 2$, then $m>10^{{10}^{9}}$.
\end{theorem}

Taking the methods of Gallot, Moree and Zudilin as in \cite{GaMoZu} further seems infeasible with current technology, as computing sufficiently many continued fractions digits is the bottleneck. However, if it would be possible to write down as many digits as one wanted, the methods in \cite{GaMoZu} lead us to believe that a lower bound for $m$ as big as $10^{{10}^{400}}$ could be achieved. Readers are encouraged to consult the excellent survey paper by Moree \cite{MoreeSurvey} for more details surrounding the mysterious Erd\H{o}s-Moser conjecture.

\medskip

 Like every other notorious conjecture in mathematics, it is also of great interest to study generalisations of the  Erd\H{o}s-Moser Conjecture.  Krzysztofek \cite{Krzysztofek} was the first to introduce an additional integer variable to right hand side of equation \eqref{eq:erdmos}.  
 
   \begin{conjecture}[generalised Erd\H{o}s-Moser, 1966]\label{conj:gEM} 
There are no integral solutions $(a,k,m)$ to the Diophantine equation 
\begin{equation}\label{eq:gEM}
S_{k}(m) = am^{k},
\end{equation}
 with $m \ge 2$, $k \ge 2$ and $a \ge 1$.
\end{conjecture}

Towards the generalised Erd\H{o}s-Moser conjecture, one can use Moser's original method to show that equation~\eqref{eq:gEM} has no integer solutions $(a,m,k)$ with 
\[
k>1, \qquad m < \max(10^{9\cdot 10^{6}}, a\cdot 10^{28}).
\] 
This result is the best known lower bound to date (see \cite{MoreeSurvey}, which improved upon bounds in \cite{MoreeGeneralised}) since we face the same obstructions as in previous works on the Erd\H{o}s-Moser Conjecture \ref{conj:ErdosMoser}; to improve on the current known lower bounds, it seems that we are required to find a completely new approach.

\medskip

The generalised Erd\H{o}s-Moser Conjecture is not the only  generalisation of  Conjecture~\ref{conj:ErdosMoser}. In another related direction, we shed light on the so called Kellner-Erd\H{o}s-Moser Conjecture \cite{Kellner}, which makes speculations on the ratios of consecutive power
sums. 

\begin{conjecture}[Kellner-Erd\H{o}s-Moser (2011)]\label{conj:KEM} 
Let $k$ and $m$ be positive integers with $m \geq  3$, then the quotient
$S_{k}(m+1)/S_{k}(m)$ is an integer if and only if  $(m, k) =(3, 1)$ or $(m, k) =(3, 3)$.
\end{conjecture}

In particular, the Kellner-Erd\H{o}s-Moser Conjecture is equivalent to 
 solving the Diophantine equation 
 \begin{equation}\label{eq:KEM}
 \alpha S_{k}(m) = m^{k},
 \end{equation}
  in the integral variables $(\alpha,k,m)$,
   where it is now apparent that we have indeed generalised Erd\H{o}s' original stipulation.

  \medskip

To contrast with equation~\eqref{eq:gEM}, equation \eqref{eq:KEM} is a bit more tame and allows one to obtain stronger results. Using a refinement of the Voronoi congruence, the notion of ``regular'' and ``irregular'' primes as well as previous techniques which relate power sums to  Bernoulli Numbers, in \cite{BaoulinaMoree2016} the authors prove the following:

\begin{theorem}[Baoulina and Moree, 2016]\label{thm:BM16}  If $\alpha$ is even or $\alpha$ has a regular prime divisor or $2 \leq \alpha \leq 1500$, then equation \ref{eq:KEM} has no solutions $(m, k)$ with $m>3$.
\end{theorem}
 
The study of  ``regular'' and ``irregular'' primes are currently of great interest to number theorists as their behaviour seems to underpin many infamous open conjectures in Number Theory.  Some well known examples are the results of Kummer relating the class number of Cyclotomic extensions with regular primes, or his proof of Fermat's Last theorem for regular prime exponents. We expect irregular primes to behave in the following manner: 
as $x \rightarrow \infty$, the number of irregular primes $p \leq x$ should be bounded above by $\delta x/(\log x)$, where $\delta$ is real constant $\delta < 1$. Conjecturally, Siegel \cite{SiegelKum} claimed that $\delta$ should be $0.3935$, which has the implication that approximately 60\% of all primes are regular. Today, Siegel's claim still remains an open problem. In fact, we don't even know if the number of regular primes is infinite!  

\medskip

 Despite the study of regular and irregular primes being highly conjectural, when pairing Theorem~\ref{thm:BM16} with expected behaviour of irregular primes, we obtain much stronger results; it then follows that, for a set of integers $\alpha$ of density 1, \eqref{eq:KEM} has no solution with $m \geq 4$.
In addition to these much stronger density-type results, Baoulina and Moree \cite{BaoulinaMoree2016} also provide explicit results for square free integers $\alpha$ whereby they give astronomical lower bounds on $k$ and $m$ for equation~\eqref{eq:KEM} having solutions.

\begin{theorem}[Baoulina and Moree, 2016] If $m \geq 4$ and $\alpha$ is square-free, then both $k$ and $m$ exceed $3.44 \cdot 10^{82}$.
\end{theorem}

\begin{table}[htbp!]
    \centering
    \begin{tabular}{ccc}
    \specialrule{.2em}{.1em}{.1em}  
        $(a, \alpha, k, m)$ & Reference(s) & Method(s) used \\
        \hline
         $(1, 1, k, m \leq 10^{10^{6}})$ & \makecell{Moser \cite{Moser1953}, \\  Moree \cite{MoreeRabbits}} & \makecell{elementary; \\ four mathemagical rabbits \\ pulled from a   tophat}\\ \hline
         $(1, 1, k \text{ odd }, m \geq 3)$ & \makecell{Moser \cite{Moser1953}, \\  Moree \cite{MoreeRabbits}} & elementary \\ \hline
         $(1, 1, k, m \leq 1.485\cdot 10^{9 321 155})$ & \makecell{Butske, Jaje \\ and Mayernick \cite{BuJaMay}} & elementary \\  \hline
         $(1, 1, k, m \leq 10^{10^{9}})$ & \makecell{Gallot, Moree \\ and Zudilin \cite{GaMoZu}} & \makecell{continued fractions \\ builds upon \cite{BestRiele76, MRU1994}}\\ \hline 
         $(a, 1, k, m < \max(10^{9\cdot 10^{6}})$ & \makecell{Moree \cite{MoreeSurvey}} & \makecell{elementary \\ builds upon \cite{MoreeGeneralised}}\\
         \hline 
         $(1, 2 \leq \alpha \leq 1500, k, m < 3)$ & \makecell{Baoulina and \\ Moree \cite{BaoulinaMoree2016}} & \makecell{elementary, \\ uses  regular primes and \\  Bernoulli numbers}\\
         \hline 
         \makecell{``For even $\alpha$ or divisible by regular \\ primes, and $m > 3$, \\ equation \eqref{eq:KEM} has  no solutions''} & \makecell{Baoulina and \\  Moree \cite{BaoulinaMoree2016}} & \makecell{elementary, \\ uses regular primes and \\  Bernoulli numbers}\\ \hline 
         \makecell{``For $\alpha$ square-free, \\ equation \eqref{eq:KEM} has  no solutions with \\ $m, k \leq 3.44\times10^{82}$''} & \makecell{Baoulina and \\ Moree \cite{BaoulinaMoree2016}} & \makecell{elementary, \\ builds upon \cite{Moser1953}}\\
         \specialrule{.2em}{.1em}{.1em} 
    \end{tabular}
    \caption{Achievements towards the resolution of equations~\eqref{eq:erdmos}, \eqref{eq:gEM} and \eqref{eq:KEM}.}
    \label{tab:EM}
\end{table}

The generalisations of equations of Erd\H{o}s-Moser type does not end here; there is a whole world of exploration unfolding from equation~\eqref{eq:erdmos}. To the reader willing to dive down this mathemagical rabbit hole, we recommend to consult \cite{Baoulina2019, MoreeSurvey, MoreeGeneralised, SonMac, SonMac1012}. Notable milestones towards the Erd\H{o}s-Moser Conjecture and its generalisations are summarised in Table~\ref{tab:EM}. We end this section with some open problems.

\subsection*{Open Problems}
\begin{itemize}
\item Develop new techniques to study \eqref{eq:erdmos}.
\item Fully resolve the infamous Erd\H{o}s-Moser Conjecture \ref{conj:ErdosMoser}. 
\item Reach the natural benchmark of $m \ge 10^{10^{7}}$ in Conjecture~\ref{conj:gEM}.
\item Density results for Conjecture~\ref{conj:gEM}.
\item Generalise and study equation~\eqref{eq:erdmos} over number fields?
\item Prove Siegel's Conjecture on the behaviour of regular/irregular primes. 
\end{itemize}

\medskip

\subsection{Sch\"affer's Conjecture} 

Sch\"affer, in 1956, taking inspiration from the speculation of Erd\H{o}s and Moser, considered the general equation 
\begin{equation} \label{eq:schaffer}
S_{k}(x) = 1^{k}+2^{k}+ \cdots + x^{k} =y^{n}. \end{equation}
Sch\"affer himself, \cite{Schaffer} gives all the (infinitely many) solutions of equation \eqref{eq:schaffer} for $(n, k)$ in the set $S := \{(1, 2), (2, 3), (3, 4), (5, 2)\}$ and proves that for all $(n, k)$ not in this set, equation \eqref{eq:table2} can only have a finite number of solutions. His proof was non-effective. However, that did not stop him from making the following conjecture.

\medskip

\begin{conjecture}[Sch\"affer] \label{conj:schaffer}
For $k \geq 1$ and $n \geq  2$ with $(k, n)$ not in the set $S$, equation \eqref{eq:table2} has only one non-trivial solution, namely
$(k, n, x, y) = (2, 2, 24, 70)$.
\end{conjecture}

\medskip

In 1979, Voorhoeve, Gy\"ory, and Tijdeman \cite{VooGyoTij} study equation \eqref{eq:table2} for $s=1$ and $r(x)$ any polynomial with integer coefficients, and prove that for $k\geq 2$ such equation has finitely many solutions.

\medskip

\begin{theorem}[Voorhoeve, Gy\"ory, and Tijdeman] Let $s=1$ and $r(x) \in \Z[x]$, for $k\geq 2$, $k\not\in \{3, 5\}$, equation \eqref{eq:table2} has finitely many integer solutions and $n$ can be effectively bounded.  
\end{theorem}

\medskip

In further results \cite{GyTiVo1980}, the same authors restrict $r(x) \in \Z$ to be an integer and present similar results for $s$ any squarefree integer. However, these are now effective. 

\medskip

\begin{theorem}[Voorhoeve, Gy\"ory, and Tijdeman] Let $s$ be a squarefree integer and $r(x) \in \Z$. For $k\geq 2$, $k\not\in \{3, 5\}$, equation \eqref{eq:table2} has finitely many integer solutions and they can be made effective.  
\end{theorem}
 
\medskip

In particular, for $s=1$ and $r=0$ they give an effective proof of Sch\"affer's Theorem whenever $k \neq 3, 5$. Their methods use a theorem of LeVeque \cite{LeVeque}, based on Siegel's method \cite{Siegel} for finding an asymptotic estimate for the frequency of solutions of a superelliptic curve.

\medskip

The effective results in \cite{GyTiVo1980} rely on the fact that for $n \geq 2$, and $r(x) \in \Z$, equation \eqref{eq:table2} gives a superelliptic equation where methods of Baker \cite{BakerHyperelliptic} based on linear forms in logarithms can be combined with LeVeque's Theorem to bound $n$, but not $x$ and $y$.

\medskip

After this initial leap, some authors worked on generalisations of the results. In \cite{Brindza1984}, Brindza considers a more general equation involving \eqref{eq:table2} and with an effective version of LeVeque's Theorem, proves an effective version of \cite{VooGyoTij}; in \cite{Brindza1990}, the same author uses similar techniques to find effective bounds for $x$ and $y$ in \cite{GyTiVo1980}; in \cite{Kano}, Kano considers more general $s$ whenever $k \not\equiv 1 \pmod{4}$; and Dilcher in \cite{Dilcher} and Urbanowicz in \cite{Urbanowicz1988}, \cite{Urbanowicz1994} and \cite{Urbanowicz1996}, introduce characters in the equation. 

\medskip

However effective, none of these results are explicit. We can find the first explicit bound for $n$ in work of Pint\'er \cite{Pinter1997}. On the other hand, Brindza \cite{Brindza1990} proves that for any given $n\neq 3$ or $4$, 
the number of integer solutions of \eqref{eq:table2}  cannot be more than $e^{7k}$. 
 Brindza and Pint\'er \cite{BriPin} further show that for $(k, n) \neq (3, 4)$, equation \eqref{eq:table2} has at most $\text{max}\{C, e^{3k}\}$ positive integer solutions when $n>2$, where $C$ is an effectively computable
absolute constant (they also give analogous results for $n=2$).

\medskip

In order to tackle Sch\"affer's conjecture, we need to completely solve equations in the family given by \eqref{eq:table2}.
Sch\"affer himself was able to confirm his conjecture for a certain finite number of pairs $(k, n)$ with $k \leq 11$, $n \geq 2$ (see \cite{Schaffer} or \cite{GyoPintSurv} for details). 
When $n=2$, Jacobson, Pint\'er and Walsh \cite{JacPintWal} confirmed Conjecture~\ref{conj:schaffer} for another handful of $k$'s, namely for all even numbers below 58 (and below 70 assuming GRH). 
Bennett, Gy\"ory and Pint\'er \cite{BenGyorPin} were the first ones to prove Sch\"affer's Conjecture for an infinite subfamily of \eqref{eq:table2}: they completely solved the problem for all $1 \leq k\leq 11$. It took all of the tools in the box to do so: from classical local and elementary methods to more powerful  modularity techniques for ternary equations, passing by linear forms in logarithms and
in elliptic logarithms, and finally appealing to computational methods such as Thue and Siegel.
Pint\'er \cite{pinter2007power}
improved on their results by showing that the conjecture still holds for $n>4$ even and $1 \leq k < 170$.

\medskip

In another direction towards tackling Sch\"affer's Conjecture, by imposing certain conditions on $x$, Hajdu \cite{Hajdu2015} recently obtained further results. Precisely, by studying the valuation of the primes 2 and 3 of the polynomial $S_{k}(x)$, Hajdu constructs upper bounds for $n$ and proves the following.

\medskip

\begin{theorem}[Hajdu]
Assume that $x \equiv 3, 4 \pmod{8}$. Then \eqref{eq:table2} has no solutions with $k = 1$
or $k$ even.
Furthermore, Sch\"affer's Conjecture holds if $x$ satisfies one of the following congruence conditions
$$x \equiv h_{i} \pmod{ m_i}, \quad \text{with } h_i \in H_i, \ \ i = 1, 2, 3, 4$$ 
where
$H_1 = \{2\}$, $H_2 = \{5, 7\}$, $H_3 = \{2, 7, 9, 14\}$, $H_4 = \{18, 22\}$,
and
$m_1 = 5$, $m_2 = 13$, $m_3 = 17$, $m_4 = 41$. 
\end{theorem}

\medskip

The theorem can be paraphrased as ``the only putative solutions to \eqref{eq:schaffer} that satisfy certain congruence conditions for $x$ are exactly the ones given by the conjectured equations''. 

\medskip

Following up on this line of thought, B\'erczes, Hajdu, Miyazaki and Pink \cite{BerHajMiyPin} managed to extend Hajdu's results and show that for all $x < 25$, all solutions of equation \eqref{eq:schaffer} with 
$n \geq 3$ are given by
$(x, k, y, n) = (1, k, 1, n), (8, 3, 6, 4)$, 
recovering Sch\"affer's observed solutions.

\medskip

The authors treat the cases where $x$ is not covered by Hajdu's Theorem \cite{Hajdu2015} via bounds from linear forms in two logarithms and
polynomial-exponential congruences. 

\medskip

Taking the ideas developed in \cite{BerHajMiyPin} further,
and combining with Baker's theory, results are obtained 
towards a \emph{generalised Sch\"affer problem}, namely
\begin{equation}\label{eq:genschaffer}
(x+1)^{k} + \cdots + (\ell x)^{k} = y^{n}  \quad k,\ell, n, x, y, s \in \Z.
\end{equation}
In \cite{BerPinSavSoy}, the authors consider $\ell=2$ in \eqref{eq:genschaffer} and for 
$2 \leq x \leq 13$, 
$k\geq 1$, 
$y\geq 2$, and 
$n\geq 3$, prove that there are no solutions to \eqref{eq:genschaffer}, as well as obtaining upper bounds
for the exponent $n$ in terms of $2$-adic and $3$-adic valuations of some functions depending on $x$ and $k$.

\medskip

In \cite{soydan2017diophantine} and \cite{BarSoy}, the authors continue studies on \eqref{eq:genschaffer}, and prove results that are in a similar vein to Sch\"affer's original problem.
The authors prove the effective finiteness of the solutions of \eqref{eq:schaffer} when taking fixed $k\geq 2$ and $\ell \geq 2$. In other words, they show that $\max\{x,y,n\} < C$, 
where $C$ is an effectively computable constant depending only on $k$ and $\ell$. Other accounts of generalisations on Sch\"affer's problem are considered in \cite{Csaba, BaPiSi}.

\medskip

Another trend on generalisations of Sch\"affer's problem was established in  \cite{bazso2012equal}, where the authors take on the endeavor of studying the polynomial
$$S_{a, b}^{k}(x):= b^k +(a+b)^k +(2a+b)^k +\ldots+(a(x-1)+b)^k.$$
Several results on the number of solution of Diophantine equations constructed upon the polynomial $S_{a, b}^{k}(x)$ are given in \cite{Bazso, BaKrLuPiRa}. 

\medskip

This closes the story of all known results on Sch\"affer's Conjecture and generalisations, for now. For a thorough early history and comprehensive study, please consult \cite{ShoreyTijdeman86}.

\begin{table}[htbp!]
    \centering
    \begin{tabular}{ccc}
    \specialrule{.2em}{.1em}{.1em}  
        $(k, n, s, r(x))$ & Reference(s) & Method(s) used \\
         \specialrule{.2em}{.1em}{.1em}  
         $(2, 2, 1, 0)$ & Watson \cite{Watson} & elliptic functions \\
         & Ljunggren \cite{Ljunggren} & algebraic proof \\ \hline
         
         $(2, 2, 1, 0)$ &  \makecell{Ma \cite{MaDiophantine}, \\ Cucurezeanu \cite{Cucurezeanu}, \\ Anglin \cite{Anglin}} & elementary \\ \hline
         
        $(2, 2, 1, 0)$ & Bennett \cite{BennettPyramid} & \ref{subsec:IntPoints}, \ref{subsec:linlogs} \\ \hline
         
         \makecell{For $(k, n) \in S$, $(s, r(x))=(1, 0)$ \eqref{eq:table2} has \\ infinitely many solutions.} & Sch\"affer \cite{Schaffer} & elementary\\ \hline
         
         \makecell{For $(k, n) \not\in S$, $(s, r(x))=(1, 0)$, \eqref{eq:table2} has \\ finitely many solutions (ineffective).} & Sch\"affer \cite{Schaffer} & \ref{subsec:Thue}, \ref{subsec:IntPoints} \\ \hline
         
         \makecell{Determination of all solutions for certain \\ small pairs $(k, n) \leq (11, 641)$; see \cite{GyoPintSurv} } & Sch\"affer \cite{Schaffer} & \ref{subsec:Thue} \\ \hline
         
         \makecell{For $k \not\in \{3, 5\}$ and $s =1$ \eqref{eq:table2} has \\ finitely many solutions (ineffective).} & \makecell{Voorhoeve, Gy\H{o}ry \\ and Tijdeman \cite{VooGyoTij}} & 
         Diophantine approximation \\ \hline
         
         \makecell{For $k \not\in \{3, 5\}$, $\square \nmid s$ and $r(x) \in \Z$, \eqref{eq:table2} has \\ finitely many solutions (effective).} & \makecell{Voorhoeve, Gy\H{o}ry \\ and Tijdeman \cite{GyTiVo1980}} & 
        \makecell{Diophantine approximation, \\ \ref{subsec:linlogs}} \\ \hline

         \makecell{Gives effective commutable  \\ bounds for $x, y$ and $n$ solutions \\ of \cite{VooGyoTij} in terms of $k$ and $s$.} & Brindza \cite{Brindza1984} &  
        \makecell{Diophantine approximation, \\ \ref{subsec:linlogs}} \\ \hline
        
         \makecell{If $n \not\in \{1, 2, 4\}$, $s=1$ and $r(x)= 0$, \\ \eqref{eq:table2} has at most $e^{7k}$ solutions} & Brindza \cite{Brindza1990} & 
         \makecell{Diophantine approximation, \\ \ref{subsec:linlogs}}\\ \hline

        $(k\leq 58, 2, 1, 0)$ & \makecell{Jacobson, Pint\'er \\ and Walsh \cite{JacPintWal}} &\ref{subsec:linlogs}, \ref{subsec:IntPoints} \\ 
        \hline
        
         $(1 \leq k \leq 11,  2\leq n, 1, 0)$ & \makecell{Bennett, Gy\H{o}ry \\ and Pint\'er \cite{BenGyorPin}} & \ref{subsec:linlogs}, \ref{subsec:Thue}, \ref{subsec:modularity} \\ \hline
         
         $(1 \leq k \leq 170 \text{ odd}, \text{even}, 1, 0)$ & Pint\'er \cite{pinter2007power} & \makecell{methods in \cite{JacPintWal} \\ with \ref{subsec:linlogs} and \ref{subsec:modularity}} \\ \hline
         
         \makecell{Solved equation \eqref{eq:table2} for a positive \\ proportion of $x$ and constructed \\ bounds on $n$ in terms of $x$ and $k$.} & Hajdu \cite{Hajdu2015} & 
         elementary
         \\ \hline
         
         \makecell{Solve \eqref{eq:table2} with $s=1$ and $r(x) =0$ \\ for $x<25$ and $n \geq 3$} & \makecell{B\'{e}rczes, Hajdu, \\ Miyazaki and Pink \cite{BerHajMiyPin}} & \makecell{methods in \cite{Hajdu2015} \\ and \ref{subsec:linlogs}} \\ 
         \specialrule{.2em}{.1em}{.1em} 
    \end{tabular}
    \caption{Progress towards Sch\"affer's Conjecture~\ref{conj:schaffer}.}
    \label{tab:Table2}
\end{table}

\subsection*{Open Problems}
\begin{itemize}
\item Solve Sch\"affer's Conjecture~\ref{conj:schaffer}.
\item Study and solve generalisations of  Sch\"affer's  Conjecture.
\end{itemize}

\section{To the power 1: integers that are sums of like powers}
A problem that may seem similar to that treated in this survey, but is actually extremely different in terms of the methodology used to tackle it, is to determine which integers 
are sums of like powers. We give a brief overview of this related problem for completeness.

\medskip

We begin by stating Fermat's infamous Two Square Theorem. 
\begin{theorem}[Fermat]
    An odd prime can be represented as the sum of two integral squares if and only if the prime is congruent 
    to 1 modulo 4.
\end{theorem}
Although the first observation (not necessarily restricted to primes) of this identity was in 1625 by Girard,
Fermat meticulously expanded upon this observation in a letter to Mersenne dated 25 December 1640: for this reason, the theorem is sometimes known as \emph{Fermat's Christmas Theorem}.

\medskip

The story then gains momentum once again in 1770, with Lagrange announcing 
his proof of the Four Square Theorem:
\begin{theorem}[Lagrange]
    Every natural number can be represented as a sum of four non-negative integral squares.
\end{theorem}
Bachet's 1621 Latin translation of Diophantus' \emph{Arithmetica} stated this theorem, however lacked a proof. Around the same time, Waring stated general observations in his Meditationes Algebraicae (pg.203--204):
\begin{displayquote}
    Omnis integer numerus vel est cubus, vel e duobus, tribus,
    4, 5, 6, 7, 8, vel novem cubus compositus: est etiam 
    quadratoquadratus; vel e duobus, tribus 
    \&c. usque ad novemdecim compositus \&sic deinceps.
\end{displayquote}
In a later edition, Waring adds (pg.349)
\begin{displayquote}
    \ldots consimilia etiam  affirmari possunt 
    (exceptis excipiendis) do oedem numero quantitatum earundem dimensionum.
\end{displayquote}
Waring's assertions/problem is more familiarly stated as
\begin{displayquote}
Can every positive integer be represented as the sum of at most $g(k)$ 
$k$th powers of positive integers, where $g(k)$ depends only on $k$ and not the number being represented? 
\end{displayquote}

\bigskip

Results so far show that every positive integer can be represented with at most 4 squares (Lagrange, 1770); at most 9 positive cubes \cite{Wieferich}, \cite{Kempner}; at most 19 fourth powers \cite{BalDesDreI}, \cite{BalDesDreII}; at most 37 positive fifth powers \cite{Chen} and at most 73 sixth powers \cite{Pillai}. 

\medskip

During the 1920's, Hardy and Littlewood introduced in a series of papers (\cite{HL2, HL3, HL4, HL5, HL6, HL8}) what we may now consider to be a more fundamental number than $g(k)$; the number $G(k)$, which is defined to be the least positive integer such that all sufficiently large integers may be represented as a sum of at most $G(k)$ $k$th powers of \textbf{positive} integers. This can be seen as a natural generalisation to Waring's Problem, or perhaps even a strengthening to Waring's Problem. 

\medskip

Since Wieferich has already shown that every positive integer can be written as a sum of  at most 9 positive cubes \cite{Wieferich}, it seems natural to study the number $G(3)$ i.e. representations of integers as sums of less than 9 cubes. 

\medskip

Landau \cite{Landau1911} showed that every sufficiently large integer is the sum of at most 8 positive cubes, however, Landau's result was not effective. Effectivity came in 1913, only two years later, where Baer \cite{Baer1913} showed that every integer greater than $2.26\times 10^{15}$ is the sum of at most 8 positive cubes. Dickson \cite{Dickson1939} refined this to the best result possible; he proves that every positive integer except  $23$ and $239$  can be written as the sum of at most 8 positive cubes. 

\medskip

As one can imagine, the problem very steeply increases with difficulty when only allowing for smaller and smaller number of positive cubes!
Linnik \cite{Linnik43} was the first, in 1943, to show that every sufficiently large integer is the sum of seven positive cubes, again, an ineffective result. A simpler, albeit still ineffective, proof was given by Watson in 1951 \cite{Watson51}. It was only until 1984 that we have effective results, the first given by McCurley \cite{McCurley84} which states that any integer greater than $\exp(\exp{13.94})$ can be written as the sum of at most 7 positive cubes.
Ramar\'e, \cite{R2007} improved this to $3.72. \times 10^{227}$. Partial results followed;
see \cite{DHL2000, BRZ1999, BElkies, Elkies}, but none settled the question once and for all.

\medskip

The question of integers represented by 7 positive cubes was finally settled in 2016 \cite{Siksek7Cubes}.

\begin{theorem}[Siksek]\label{thm:Siksekcubes}
Every positive integer except
\[
15, 22, 23, 50, 114, 167, 175, 186, 212, 231, 238, 239, 303, 364, 420, 428, 454
\]
is the sum of seven cubes.
\end{theorem}

Currently, we know that $4 \leq G(3) \leq 7$. It is conjectured that $G(3)=4$, however, we are yet to determine whether every sufficiently large integer can be represented as the sum of 6 positive cubes, let alone specific results along the lines of Theorem~\ref{thm:Siksekcubes}. In Table~\ref{tab:Table6} we summarise the known bounds on $G(k)$. An excellent survey article on Waring's Problem and its generalisations is \cite{VW2000} and the references therein, where there are more results on the function $G(k)$.

\begin{table}[htbp!]
    \centering
    \begin{tabular}{ccccc}
    \specialrule{.2em}{.1em}{.1em}  
       &   & Known Bounds on $G(k)$ & & \\
         \specialrule{.2em}{.1em}{.1em} 
      $1 = G(1) = 1$ & 
      $6 \leq G(5) \leq 17 $&
      $13 \leq G(9) \leq 50 $& 
      $14 \leq G(13) \leq 84 $ & 
      $18 \leq G(17) \leq 117 $ \\
      
      $4 =G(2) = 4$&
      $9 \leq G(6) \leq 24 $ &
      $12\leq G(10) \leq 59 $ &
      $15 \leq G(14) \leq 92 $ & 
      $27 \leq G(18) \leq 125 $ \\
      
      $4 \leq G(3) \leq 7$&
      $8 \leq G(7) \leq 33 $ & 
      $12 \leq G(11) \leq 67 $& 
      $16 \leq G(15) \leq 100 $&  
      $20 \leq G(19) \leq 134 $ \\
      
      $16 =G(4) = 16$&
      $32 \leq G(8) \leq 42 $& 
      $16 \leq G(12) \leq 76 $& 
      $64 \leq G(16) \leq 109$& 
      $25 \leq G(20) \leq 142 $ \\
         \specialrule{.2em}{.1em}{.1em} 
    \end{tabular}
    \caption{Known bounds on $G(k)$.}
    \label{tab:Table6}
\end{table}

\bigskip

We now shift gears slightly, to talk about \emph{explicit} representations of certain integers as the sum of a specified number of $k$th powers - we are able to speak of such things as we no longer specify that the $k$th powers need to be positive. In this direction of study, Lagrange's results were extended by Legendre in 1797-1798, who considered integers that are sums of three squares:
\begin{theorem}[Legendre]
    A natural number can be represented as the sum of three squares if and only if it is not of the form
    $4^{a}(8b+7)$ for some integers $a$ and $b$.
\end{theorem}

\medskip

In determining integers that are sums of a small number of (integer) cubes, the problem becomes significantly more difficult. We begin gently with the problem of expressing integers as a sum of two integer cubes, in other words,
let $k$ be an arbitrary integer, then, is it possible to write 
\[ 
k = x^3 + y^3, \qquad \text{ for some } x,y \in \Z?
\]
This problem is actually not too hard to solve, since we have $x^3+y^3 = (x+y)(x^2-xy+y^2)=rs$, with $r=x+y$, and from this we obtain $s=3x^2 -3rx+r^2$, which has integer solutions if and only if $12s-3r^2$ is a square, say $t^2$, and in particular the solutions are $x=(3r+t)/6,\ y=(3r-t)/6$.

\medskip

It is worth noticing that, if instead of integer $x$ and $y$ we allow for rational numbers, the problem of solving $k=x^3+y^3$ for a given $k$ is transforms to that of finding rational points on an elliptic curve. 
The curve $x^3+y^3=k$ has Weierstrass form $Y^2=X^3-432k^2$, it has nontrivial torsion if and only if $k$ is a cube or twice a cube, and otherwise we have a solution if and only if it has positive rank. There are some partial results towards this, which we list below.
\begin{itemize}
    \item In \cite{ABS}, it is proven that a positive proportion of $k$ yield a rank zero elliptic curve and a positive proportion of $k$ yield a positive rank elliptic curve.
    \item In \cite{satge86, elkies94, DV17, kriz20} there are results concerning the case where $k$ is a prime number, namely if $k\equiv 2,5 \pmod 9$ then it yields a rank zero elliptic curve, while if $k \equiv 4,7,8 \pmod 9$ it yields a rank 1 elliptic curve. For $k \equiv 1 \pmod 9$, work of \cite{RVZ95} shows that the underlying elliptic curve may have rank either 0 or 2.
    \item In \cite{KL20}, the cases where $k$ is twice a prime $p$ or twice the square of a prime $p$ satisfying $p\equiv 2,5 \pmod 9$ are dealt with by using Iwasawa theory. For a more informative but less technical reading on the subject, it is worth consulting \cite{Kezuka22} (in Japanese).
\end{itemize}

\medskip
 
A well-known result concerns the case of the sum of five or more cubes: indeed every integer has infinitely many representations as the sum of five cubes, since the following identity holds:
$$ 6m=(m+1)^3 +(m-1)^3 -2m^3 $$
and 
if we write $k=6n+r$, we can apply this identity to $m=m(n) = (k-(6n+r)^3)/6$, obtaining:
$$ k=(6n+r)^3+(m(n)+1)^3+(m(n)-1)^3-(m(n))^3-(m(n))^3. $$

\medskip

In \cite{Dem66} it is shown, with a similar approach, that every integer $k \not\equiv \pm 4 \pmod 9$ can be represented as a sum of four cubes in infinitely many ways, and it is conjectured (see \cite{Sier60}) that in fact every integer $k$ can be written as a sum of four cubes in infinitely many ways.

\medskip

Finally, let us consider the problem of whether a given integer can be written as a sum of three cubes. First of all, a simple congruence condition modulo $9$ shows that if $k \equiv \pm 4 \pmod 9$ then it is not the sum of three cubes. Moreover, it is easy to show that if $k$ is a cube or twice a cube then again there are infinitely many ways to write it as a sum of three cubes -- it is sufficient to multiply by a suitable cube one of the following identities:
\begin{align*}
0 &=n^3-n^3+0^3;\\
1 &= (9n^4)^3 + (3n-9n^4)^3 +(1-9n^3)^3;\\
2 &= (1+6n^3)^3 + (1-6n^3)^3 +(-6n^2)^3.
\end{align*}

\medskip

It remains to understand whether integers $k\not\equiv \pm 4 \pmod 9$ and not being a cube or twice a cube can be written as the sum of three cubes, and if so, whether there are infinitely many ways of doing so. In 1956 \cite{Mordell56}, Mordell wrote: \textit{I do not know anything about the integer solutions of $X^3+Y^3+Z^3=3$ beyond the existence of $(1,1,1),\ (4,4,-5)$}, and indeed this was such a hard problem that only over 60 years later a less trivial solution was found \cite{BS21} by Booker and Sutherland, namely:
$$ 569936821221962380720^3+(-569936821113563493509)^3+(-472715493453327032)^3 =3. $$

\medskip

They also treat rigorously and computationally the challenge of writing, more generally, an integer $k$ as the sum of three cubes, other than the cases that are already understood. Although it may seem an innocuous computational task, it is incredibly hard to compute numbers with so many digits by a na\"ive brute force algorithm. Indeed, already in 1954, Miller and Woollett \cite{MW54} programmed the EDSAC (Electronic Delay Storage Automatic Calculator) in Cambridge to conduct the search for solutions to Mordell's equation with $0 \leq k \leq 100$ and $|Z|\leq |Y| \leq |X| \leq 3164$, and succeeded in finding solutions for $k\not\equiv \pm 4 \pmod 9 $ in this range, with the exceptions of 
$$ k = 30,33,39,42,52,74,75,84,87. $$

\medskip

Later, in 1964, a more ambitious and sophisticated search \cite{Gardiner64} involved all $0 \leq k \leq 1000$ and $|Z|\leq |Y| \leq |X| \leq 2^{16}$, but was only able to produce one extra solution, for $k=87$.

\medskip

In 1992, Heath-Brown \cite{HB92} finally made a more precise conjecture, namely that every integer $k\not\equiv \pm 4 \pmod 9$ has infinitely many writings as sum of three cubes. He also conjectured an asymptotic formula for this number of solutions, as follows: if $N(B)$ is the number of solutions in the box $\max \{|X|,|Y|,|Z| \}\leq B$, then
$$ N(B) \sim \dfrac{1}{9} \dfrac{\Gamma(1/3)^2}{\Gamma(2/3)}\left(\prod_{p \text{ prime}} \sigma_p \right) \log B, $$
where $\sigma_p $ is defined, for example, in \cite{SiksekSurvey}.

\medskip

Due to the above partial results and significant progress towards the conjectures, there has been strong motivation towards finding more optimised algorithms to search for solutions to Mordell's equation (also where $3$ is replaced by an integer $k$). Following \cite{SiksekSurvey}, we illustrate here the most significant algorithm currently known, which was used, among others, by Booker and Sutherland to find the result mentioned above, as well as to write every number below $100$ (with the usual restrictions) as a sum of three cubes. This algorithm is based on the ideas and work of \cite{HB,HBetal,Beck,Booker,BS21} and it is somewhat similar to the factorisation algorithm mentioned above for the case of two cubes.

\medskip

Fix $k$, suppose there is a solution $(X,Y,Z)$ to $X^3+Y^3+Z^3=k$ and let $d=X+Y$. Then, taking $X=d-Y$ we obtain a quadratic equation for $Y$ in terms of $d,Z$, such that its discriminant is a square. Therefore, we have two constraints on $Z$, namely:
\begin{align*}
    Z^3 & \equiv k \pmod d,\\
    3d(4(k&-Z^3)-d^3) \text{ is a square.}
\end{align*}

\medskip

At this point, Booker \cite{Booker} computes cube roots of $k$ modulo $d$ for $d$ varying over a large range of values, ordered by prime factorisation (rather than running the algorithm on consecutive $d$'s). This way, it is possible to find solutions in the box $\max\{|X|,|Y|,|Z|\} \leq B$ in $O(B \log \log B \log \log \log B)$ arithmetic operations and table look-ups.

\medskip

To end this section, we include a table with the solutions to $X^3+Y^3+Z^3=k$ for $k \in \{30,33,39,42,52,74,75,84,87\}$. For references about who the first authors to find such solutions were, see \cite[Table 1]{SiksekSurvey}. We list them here in order of ``discovery''.

\begin{table}[htbp!]
    \centering
    \begin{tabular}{cc}
    \specialrule{.2em}{.1em}{.1em} 
        $k$ & $X,Y,Z$ \\
        \hline
        87 & 4271, -4126, -1972 \\ 
        39 & 134476, 117367, -159380\\
        84 & 41639611, -41531726, -8241191\\
        75 & 4381159, 435203083, -435203231\\
        30 & 2220422932, -283059965, -2218888517\\
        52 & 60702901317, 23961292454, -61922712865\\
        74 & -284650292555885, 66229832190556, 283450105697727\\
        33 & 8866128975287528, -8778405442862239, -2736111468807040\\
        42 & -80538738812075974, 80435758145817515, 12602123297335631\\
        \specialrule{.2em}{.1em}{.1em}
    \end{tabular}
    \caption{Last remaining numbers to be written as a sum of three cubes}
    \label{tab:sums3cubes}
\end{table}

\subsection*{Open Problems}
\begin{itemize}
\item Partial, or even full progress on Waring's Problem.
\item Can every sufficiently large integer be represented as the sum of 6 positive cubes?
\item Can every sufficiently large integer be represented as the sum of 5 positive cubes?
\item Explicit determination of integers that are a small number of cubes.
\end{itemize}

\section*{Acknowledgements}

This collaboration grew from the WiNE-4 workshop that was hosted at Universiteit Utrecht in August 2022. The authors would like to express their sincere gratitude to the organisers of the workshop; Ramla Abdellatif,
Valentijn Karemaker,
Ariane M\'ezard and Nirvana Coppola,
 for creating a welcoming and hospitable environment that resulted in a fruitful and thoroughly enjoyable meeting.
The authors would like to thank the anonymous referee for carefully reading the paper and suggesting several improvements. The authors would like to extend their thanks to  Pieter Moree, and \'Akos Pint\'er for pointing out key references.

N. Coppola is supported by the NWO Vidi grant No. 639.032.613, New Diophantine Directions. 
M. Khawaja is supported by an EPSRC studentship from the University of Sheffield (EPSRC grant no. EP/T517835/1). 
\"{O}. \"{U}lkem is supported by T\"{U}BITAK project no. 119F405.

\bibliographystyle{siam}
\bibliography{ref.bib}

\begin{thebibliography}{100}

\bibitem{ABS}
{\sc L.~Alp\"{o}ge, M.~Bhargava, and A.~Shnidman}, {\em Integers expressible as
  the sum of two rational cubes (with an appendix by {A}. {B}urungale and {C}.
  {S}kinner)}, \href{https://arxiv.org/abs/2210.10730}{arXiv: 2210.10730},
  (2022).

\bibitem{Anglin}
{\sc W.~S. Anglin}, {\em The square pyramid puzzle}, Amer. Math. Monthly,
  \textbf{97} (1990), pp.~120--124.

\bibitem{Garcia2019}
{\sc A.~Arg\'{a}ez-Garc\'{i}a}, {\em On perfect powers that are sums of cubes
  of a five term arithmetic progression}, J. Number Theory, \textbf{201}
  (2019), pp.~460--472.

\bibitem{GarciaPatel2019}
{\sc A.~Arg\'{a}ez-Garc\'{i}a and V.~Patel}, {\em On perfect powers that are
  sums of cubes of a three term arithmetic progression}, J. Comb. Number
  Theory, \textbf{10}(3) (2019), pp.~147--160.

\bibitem{GarciaPatel2020}
\leavevmode\vrule height 2pt depth -1.6pt width 23pt, {\em On perfect powers
  that are sums of cubes of a seven term arithmetic progression}, J. Number
  Theory, \textbf{214} (2020), pp.~440--451.

\bibitem{Baer1913}
{\sc W.~S. Baer}, {\em \"{U}ber die {Z}erlegung der ganzen {Z}ahlen in sieben
  {K}uben}, Math. Ann., \textbf{74} (1913), pp.~511--514.

\bibitem{Baker66}
{\sc A.~Baker}, {\em Linear forms in logarithms of algebraic numbers},
  Mathematika, \textbf{13} (1966), pp.~204--216.

\bibitem{BakerHyperelliptic}
\leavevmode\vrule height 2pt depth -1.6pt width 23pt, {\em Bounds for the
  solutions of the hyperelliptic equation}, Proc. Cambridge Philos. Soc.,
  \textbf{65} (1969), pp.~439--444.

\bibitem{baker90}
\leavevmode\vrule height 2pt depth -1.6pt width 23pt, {\em Transcendental
  number theory}, Cambridge Mathematical Library, Cambridge University Press,
  Cambridge, second edition~ed., 1990.

\bibitem{BalDesDreI}
{\sc R.~Balasubramanian, J.-M. Deshouillers, and F.~Dress}, {\em Probl\`{e}me
  de {W}aring pour les bicarr\'{e}s. {I}. {S}ch\'{e}ma de la solution}, C. R.
  Acad. Sci. Paris S\'{e}r. I Math., \textbf{303} (1986), pp.~85--88.

\bibitem{BalDesDreII}
\leavevmode\vrule height 2pt depth -1.6pt width 23pt, {\em Probl\`eme de
  {W}aring pour les bicarr\'{e}s. {II}. {R}\'{e}sultats auxiliaires pour le
  th\'{e}or\`{e}me asymptotique}, C. R. Acad. Sci. Paris S\'{e}r. I Math.,
  \textbf{303} (1986), pp.~161--163.

\bibitem{Baoulina2019}
{\sc I.~N. Baoulina}, {\em On the unsolvability of certain equations of
  {E}rd{\H o}s-{M}oser type}, Elem. Math., \textbf{74} (2019), pp.~1--9.

\bibitem{BaoulinaMoree2016}
{\sc I.~N. Baoulina and P.~Moree}, {\em Forbidden integer ratios of consecutive
  power sums}, in From arithmetic to zeta-functions, Springer, [Cham], 2016,
  pp.~1--30.

\bibitem{BarSoy}
{\sc D.~Bartoli and G.~Soydan}, {\em The {D}iophantine equation
  {$(x+1)^k+(x+2)^k +\dots+(\ell x)^k= y^n$} revisited}, Publ. Math. Debrecen,
  \textbf{96} (2020), pp.~111--120.

\bibitem{Bazso}
{\sc A.~Bazs\'o}, {\em Effective results for polynomial values of (alternating)
  power sums of arithmetic progressions}, preprint,  (2020).

\bibitem{BazBerGyoPint}
{\sc A.~Bazs\'{o}, A.~B\'{e}rczes, K.~Gy{\H o}ry, and A.~Pint\'{e}r}, {\em On
  the resolution of equations {$Ax^n-By^n=C$} in integers {$x,y$} and
  {$n\geq3$}. {II}}, Publ. Math. Debrecen, \textbf{76} (2010), pp.~227--250.

\bibitem{bazso2012equal}
{\sc A.~Bazs{\'o}, D.~Kreso, F.~Luca, and {\'A}.~Pint{\'e}r}, {\em On equal
  values of power sums of arithmetic progressions}, Glasnik matemati{\v{c}}ki,
  \textbf{47} (2012), pp.~253--263.

\bibitem{BaKrLuPiRa}
{\sc A.~Bazs\'o, D.~Kreso, F.~Luca, A.~Pint\'er, and {\relax Cs}.~Rakaczki},
  {\em On equal values of products and power sums of consecutive elements in an
  arithmetic progression}, preprint,  (2023).

\bibitem{BaPiSi}
{\sc A.~Bazs\'{o}, A.~Pint\'{e}r, and H.~M. Srivastava}, {\em A refinement of
  {F}aulhaber's theorem concerning sums of powers of natural numbers}, Appl.
  Math. Lett., \textbf{25} (2012), pp.~486--489.

\bibitem{Beck}
{\sc M.~Beck, E.~Pine, W.~Tarrant, and K.~Yarbrough~Jensen}, {\em New integer
  representations as the sum of three cubes}, Math. Comp., \textbf{76} (2007),
  pp.~1683--1690.

\bibitem{BennettPyramid}
{\sc M.~A. Bennett}, {\em Lucas' square pyramid problem revisited}, Acta
  Arith., \textbf{105} (2002), pp.~341--347.

\bibitem{BenDahmen}
{\sc M.~A. Bennett and S.~R. Dahmen}, {\em Klein forms and the generalized
  superelliptic equation}, Ann. of Math. (2), \textbf{177} (2013),
  pp.~171--239.

\bibitem{BGMP2006}
{\sc M.~A. Bennett, K.~Gy{\H o}ry, M.~Mignotte, and A.~Pint\'{e}r}, {\em
  Binomial {T}hue equations and polynomial powers}, Compos. Math., 142 (2006),
  pp.~1103--1121.

\bibitem{BenGyorPin}
{\sc M.~A. Bennett, K.~Gy{\H o}ry, and A.~Pint\'{e}r}, {\em On the
  {D}iophantine equation {$1^k+2^k+\dots+x^k=y^n$}}, Compos. Math.,
  \textbf{140} (2004), pp.~1417--1431.

\bibitem{BennettKoutsianas}
{\sc M.~A. Bennett and A.~Koutsianas}, {\em The equation
  {$(x-d)^5+x^5+(x+d)^5=y^n$}}, Acta Arith., \textbf{198} (2021), pp.~387--399.

\bibitem{Bennett2016}
{\sc M.~A. Bennett, V.~Patel, and S.~Siksek}, {\em {S}uperelliptic equations
  arising from sums of consecutive powers}, Acta Arith., \textbf{172} (2016),
  pp.~377--393.

\bibitem{bennett2017perfect}
\leavevmode\vrule height 2pt depth -1.6pt width 23pt, {\em Perfect powers that
  are sums of consecutive cubes}, Mathematika, \textbf{63} (2017),
  pp.~230--249.

\bibitem{BM2}
{\sc M.~A. Bennett and S.~Siksek}, {\em Differences between perfect powers :
  prime power gaps}, Algebra Number Theory, \emph{to appear} (2023).

\bibitem{BM1}
\leavevmode\vrule height 2pt depth -1.6pt width 23pt, {\em Differences between
  perfect powers: the {L}ebesgue-{N}agell equation}, Trans. Amer. Math. Soc.,
  \textbf{376} (2023), pp.~335--370.

\bibitem{BenSki}
{\sc M.~A. Bennett and C.~M. Skinner}, {\em Ternary {D}iophantine equations via
  {G}alois representations and modular forms}, Canad. J. Math., \textbf{56}
  (2004), pp.~23--54.

\bibitem{BenVatYaz}
{\sc M.~A. Bennett, V.~Vatsal, and S.~Yazdani}, {\em Ternary {D}iophantine
  equations of signature {$(p,p,3)$}}, Compos. Math., \textbf{140} (2004),
  pp.~1399--1416.

\bibitem{BerHajMiyPin}
{\sc A.~B\'{e}rczes, L.~Hajdu, T.~Miyazaki, and I.~Pink}, {\em On the equation
  {$1^k+2^k+\dots+x^k=y^n$} for fixed {$x$}}, J. Number Theory, \textbf{163}
  (2016), pp.~43--60.

\bibitem{BerPinSavSoy}
{\sc A.~B\'{e}rczes, I.~Pink, G.~Sava\c{s}, and G.~Soydan}, {\em On the
  {D}iophantine equation {$(x+1)^k+(x+2)^k+\dots+(2x)^k=y^n$}}, J. Number
  Theory, \textbf{183} (2018), pp.~326--351.

\bibitem{BRZ1999}
{\sc F.~Bertault, O.~Ramar\'{e}, and P.~Zimmermann}, {\em On sums of seven
  cubes}, Math. Comp., 68 (1999), pp.~1303--1310.

\bibitem{BestRiele76}
{\sc M.~R. Best and H.~J.~J. te~Riele}, {\em On a conjecture of {E}rd{\H o}s
  concerning sums of powers of integers}, Stichting Mathematisch Centrum.
  Numerieke Wiskunde,  (1976), p.~Report NW 23/76.

\bibitem{BHV}
{\sc {\relax Yu}.~Bilu, G.~Hanrot, and P.~M. Voutier}, {\em Existence of
  primitive divisors of {L}ucas and {L}ehmer numbers (with an appendix by {M}.
  {M}ignotte)}, J. Reine Angew. Math, \textbf{539} (2001), pp.~75--122.

\bibitem{BElkies}
{\sc K.~D. Boklan and N.~D. Elkies}, {\em Every multiple of $4$ except $212,
  364, 420$, and $428$ is the sum of seven cubes},
  \href{https://arxiv.org/abs/0903.4503}{arXiv: 0903.4503},  (2009).

\bibitem{Booker}
{\sc A.~R. Booker}, {\em Cracking the problem with {$33$}}, Res. Number Theory,
  \textbf{5} (2019), pp.~1--6.

\bibitem{BS21}
{\sc A.~R. Booker and A.~V. Sutherland}, {\em On a question of {M}ordell},
  Proc. Natl. Acad. Sci. U.S.A., \textbf{118} (2021), p.~e2022377118.

\bibitem{Magma}
{\sc W.~Bosma, J.~Cannon, and C.~Playoust}, {\em The {M}agma algebra system.
  {I}. {T}he user language}, J. Symbolic Comput., \textbf{24} (1997),
  pp.~235--265.
\newblock Computational algebra and number theory (London, 1993).

\bibitem{box2021elliptic}
{\sc J.~Box}, {\em Elliptic curves over totally real quartic fields not
  containing {$\sqrt{5}$} are modular}, Trans. Amer. Math. Soc., \textbf{375}
  (2022), pp.~3129--3172.

\bibitem{BreStrTza}
{\sc A.~Bremner, R.~J. Stroeker, and N.~Tzanakis}, {\em On sums of consecutive
  squares}, J. Number Theory, \textbf{62} (1997), pp.~39--70.

\bibitem{BCDT2001}
{\sc C.~Breuil, B.~Conrad, F.~Diamond, and R.~Taylor}, {\em On the {M}odularity
  of {E}lliptic {C}urves over {$\mathbb{Q}$}: {W}ild 3-{A}dic {E}xercises}, J.
  Amer. Math. Soc, \textbf{14} (2001), pp.~843--939.

\bibitem{Brindza1984}
{\sc B.~Brindza}, {\em On some generalizations of the {D}iophantine equation
  {$1^k+2^k+\cdots +x^k=y^z$}}, Acta Arith., \textbf{44} (1984), pp.~99--107.

\bibitem{Brindza1990}
\leavevmode\vrule height 2pt depth -1.6pt width 23pt, {\em Power values of sums
  {$1^k+2^k+\cdots+x^k$}}, in Number theory, {V}ol. {II} ({B}udapest, 1987),
  vol.~\textbf{51} of Colloq. Math. Soc. J\'{a}nos Bolyai, North-Holland,
  Amsterdam, 1990, pp.~595--611.

\bibitem{BriPin}
{\sc B.~Brindza and A.~Pint\'{e}r}, {\em On the number of solutions of the
  equation {$1^k+2^k+\cdots+(x-1)^k=y^z$}}, Publ. Math. Debrecen, \textbf{56}
  (2000), pp.~271--277.

\bibitem{BugMigSik}
{\sc Y.~Bugeaud, M.~Mignotte, and S.~Siksek}, {\em A multi-{F}rey approach to
  some multi-parameter families of {D}iophantine equations}, Canad. J. Math.,
  \textbf{60} (2008), pp.~491--519.

\bibitem{BuJaMay}
{\sc W.~Butske, L.~M. Jaje, and D.~R. Mayernik}, {\em On the equation
  {$\sum_{p\mid N}(1/p)+(1/N)=1$}, pseudoperfect numbers, and perfectly
  weighted graphs}, Math. Comp., \textbf{69} (2000), pp.~407--420.

\bibitem{CN23}
{\sc A.~Caraiani and J.~Newton}, {\em On the modularity of elliptic curves over
  imaginary quadratic fields}, \href{https://arxiv.org/abs/2301.10509}{arXiv:
  2301.10509},  (2023).

\bibitem{Cassels}
{\sc J.~W.~S. Cassels}, {\em A {D}iophantine equation}, Glasgow Math. J.,
  \textbf{27} (1985), pp.~11–--18.

\bibitem{Chen}
{\sc J.~Chen}, {\em Waring's problem for {$g(5)=37$}}, Sci. Sinica, \textbf{13}
  (1964), pp.~1547--1568.

\bibitem{cohenVI}
{\sc H.~Cohen}, {\em Number Theory: Volume I: Tools and Diophantine Equations},
  vol.~\textbf{239} of Grad. Texts in Math., Springer New York, 2007.

\bibitem{Cohen}
{\sc H.~Cohen}, {\em Number theory: Volume II: Analytic and modern tools},
  vol.~\textbf{240} of Grad. Texts in Math., Springer, New York, 2007.

\bibitem{cohn_1996}
{\sc J.~H.~E. Cohn}, {\em Perfect {P}ell {P}owers}, Glasg. Math. J.,
  \textbf{38} (1996), pp.~19--–20.

\bibitem{CIKPU}
{\sc N.~Coppola, M.~Curc\'o-Iranzo, M.~Khawaja, V.~Patel, and O.~\"{U}lkem},
  {\em On perfect powers that are sums of cubes of a nine term arithmetic
  progression}, \href{https://arxiv.org/abs/2307.01815}{arXiv: 2307.01815},
  (2023).

\bibitem{SengunSiksek18}
{\sc M.~H. \c{S}eng\"{u}n and S.~Siksek}, {\em On the asymptotic {F}ermat's
  last theorem over number fields}, Comment. Math. Helv., 93 (2018),
  pp.~359--375.

\bibitem{Turcas18}
{\sc G.~C. \c{T}urca\c{s}}, {\em On {F}ermat's equation over some quadratic
  imaginary number fields}, Res. Number Theory, 4 (2018), pp.~Paper No. 24, 16.

\bibitem{Turcas20}
\leavevmode\vrule height 2pt depth -1.6pt width 23pt, {\em On {S}erre's
  modularity conjecture and {F}ermat's equation over quadratic imaginary fields
  of class number one}, J. Number Theory, 209 (2020), pp.~516--530.

\bibitem{Cucurezeanu}
{\sc I.~Cucurezeanu}, {\em An elementary solution of {L}ucas' problem}, J.
  Number Theory, \textbf{44} (1993), pp.~9--12.

\bibitem{DDKT}
{\sc P.~Das, P.~K. Dey, A.~Koutsianas, and N.~Tzanakis}, {\em Perfect powers in
  sum of three fifth powers}, J. Number Theory, \textbf{236} (2022),
  pp.~443--462.

\bibitem{DV17}
{\sc S.~Dasgupta and J.~Voight}, {\em Sylvester's problem and mock {H}eegner
  points}, Proc. Amer. Math. Soc., \textbf{146} (2018), pp.~3257--3273.

\bibitem{Deconinck15}
{\sc H.~Deconinck}, {\em On the generalized {F}ermat equation over totally real
  fields}, Acta Arith., \textbf{173} (2016), pp.~225--237.

\bibitem{Centina}
{\sc A.~Del~Centina}, {\em Unpublished manuscripts of {S}ophie {G}ermain and a
  revaluation of her work on {F}ermat’s {L}ast {T}heorem}, Arch. Hist. Exact
  Sci., \textbf{62} (2008), pp.~349--392.

\bibitem{Dem66}
{\sc V.~A. Dem'janenko}, {\em Sums of four cubes}, Izv. Vys\v{s}. U\v{c}ebn.
  Zaved. Matematika,  (1966), pp.~64--69.

\bibitem{cubicmod}
{\sc M.~Derickx, F.~Najman, and S.~Siksek}, {\em Elliptic curves over totally
  real cubic fields are modular}, Algebra Number Theory, \textbf{14} (2020),
  pp.~1791--–1800.

\bibitem{DHL2000}
{\sc J.-M. Deshouillers, F.~Hennecart, and B.~Landreau}, {\em
  {$7\,373\,170\,279\,850$}}, Math. Comp., \textbf{69} (2000), pp.~421--439.

\bibitem{Dickson1939}
{\sc L.~E. Dickson}, {\em All integers except {$23$} and {$239$} are sums of
  eight cubes}, Bull. Amer. Math. Soc., \textbf{45} (1939), pp.~588--591.

\bibitem{Dickson}
\leavevmode\vrule height 2pt depth -1.6pt width 23pt, {\em History of the
  {T}heory of {N}umbers, vol. 2}, Diophantine Analysis, Chelsea, New York,
  1971.

\bibitem{Dilcher}
{\sc K.~Dilcher}, {\em On a {D}iophantine equation involving quadratic
  characters}, Compos. Math., \textbf{57} (1986), pp.~383--403.

\bibitem{EdisThesis}
{\sc S.~Edis}, {\em On Arithmetic Progressions and Perfect Powers}, PhD thesis,
  University of Sheffield, 2019.

\bibitem{elkies94}
{\sc N.~D. Elkies}, {\em Heegner point computations}, in Algorithmic Number
  Theory, Springer, Berlin, Heidelberg, 1994, pp.~122--133.

\bibitem{Elkies}
\leavevmode\vrule height 2pt depth -1.6pt width 23pt, {\em Every even number
  greater than 454 is the sum of seven cubes},
  \href{https://arxiv.org/abs/1009.3983}{arXiv: 1009.3983},  (2010).

\bibitem{EulerAlg}
{\sc L.~Euler}, {\em Vollst\"{a}ndige {A}nleitung zur {A}lgebra},
  Reclam-Verlag, Stuttgart, 1959.

\bibitem{Faltings}
{\sc G.~Faltings}, {\em Diophantine approximation on abelian varieties}, Ann.
  of Math., \textbf{133} (1991), pp.~549--576.

\bibitem{freitas2014elliptic}
{\sc N.~Freitas, B.~V.~L. Hung, and S.~Siksek}, {\em Elliptic curves over real
  quadratic fields are modular}, Invent. Math., \textbf{201} (2015),
  pp.~159--206.

\bibitem{FREITAS2020106964}
{\sc N.~Freitas, A.~Kraus, and S.~Siksek}, {\em Class field theory,
  {D}iophantine analysis and the asymptotic {F}ermat's last theorem}, Adv.
  Math., \textbf{363} (2020), pp.~106964, 37.

\bibitem{Freitas}
{\sc N.~Freitas and S.~Siksek}, {\em Fermat's {L}ast {T}heorem over some small
  real quadratic fields}, Algebra Number Theory, \textbf{9} (2014),
  pp.~875--895.

\bibitem{freitas_siksek_2015}
{\sc N.~Freitas and S.~Siksek}, {\em The asymptotic {F}ermat’s {L}ast
  {T}heorem for five-sixths of real quadratic fields}, Compos. Math.,
  \textbf{151} (2015), pp.~1395--–1415.

\bibitem{GaMoZu}
{\sc Y.~Gallot, P.~Moree, and W.~Zudilin}, {\em The {E}rd{\H o}s-{M}oser
  equation {$1^k+2^k+\dots+(m-1)^k=m^k$} revisited using continued fractions},
  Math. Comp., \textbf{80} (2011), pp.~1221--1237.

\bibitem{Gardiner64}
{\sc V.~Gardiner, R.~B. Lazarus, and P.~R. Stein}, {\em Solutions of the
  {D}iophantine equation {$x^3+y^3=z^3-d$}}, Math. Comp., \textbf{18} (1964),
  pp.~408--413.

\bibitem{GouWang}
{\sc S.~Gou and T.~Wang}, {\em The {D}iophantine equation
  {$x^2+2a\cdot17b=y^n$}}, Czechoslovak Math. J., \textbf{62} (2012),
  pp.~645--654.

\bibitem{GyoPintSurv}
{\sc K.~Gy{\H o}ry and A.~Pint\'{e}r}, {\em On the equation
  {$1^k+2^k+\dots+x^k=y^n$}}, Publ. Math. Debrecen, \textbf{62} (2003),
  pp.~403--414.

\bibitem{GyTiVo1980}
{\sc K.~Gy{\H o}ry, R.~Tijdeman, and M.~Voorhoeve}, {\em On the equation
  {$1^{k}+2^{k}+\cdots +x^{k}=y^{z}$}}, Acta Arith., \textbf{37} (1980),
  pp.~233--240.

\bibitem{Hajdu2015}
{\sc L.~Hajdu}, {\em On a conjecture of {S}ch\"{a}ffer concerning the equation
  {$1^k+\dots+x^k=y^n$}}, J. Number Theory, \textbf{155} (2015), pp.~129--138.

\bibitem{HL2}
{\sc G.~H. Hardy and J.~E. Littlewood}, {\em Some problems of ``partitio
  numerorum'': {II}. {P}roof that every large number is the sum of at most 21
  biquadrates}, Math. Z., \textbf{9} (1921), pp.~14--27.

\bibitem{HL4}
\leavevmode\vrule height 2pt depth -1.6pt width 23pt, {\em Some problems of
  `{P}artitio {N}umerorum': {IV}. {T}he singular series in {W}aring's {P}roblem
  and the value of the number {$G(k)$}}, Math. Z., \textbf{12} (1922),
  pp.~161--188.

\bibitem{HL3}
\leavevmode\vrule height 2pt depth -1.6pt width 23pt, {\em Some problems of
  `{P}artitio numerorum'; {III}: {O}n the expression of a number as a sum of
  primes}, Acta Math., \textbf{44} (1923), pp.~1--70.

\bibitem{HL5}
\leavevmode\vrule height 2pt depth -1.6pt width 23pt, {\em Some {P}roblems of
  '{P}artitio {N}umerorum'({V}): {A} {F}urther {C}ontribution to the {S}tudy of
  {G}oldbach's {P}roblem}, Proc. London Math. Soc. (2), \textbf{22} (1924),
  pp.~46--56.

\bibitem{HL6}
\leavevmode\vrule height 2pt depth -1.6pt width 23pt, {\em Some problems of
  `{P}artitio numerorum' ({VI}): {F}urther researches in {W}aring's {P}roblem},
  Math. Z., \textbf{23} (1925), pp.~1--37.

\bibitem{HL8}
\leavevmode\vrule height 2pt depth -1.6pt width 23pt, {\em Some {P}roblems of
  '{P}artitio {N}umerorum' ({VIII}): {T}he {N}umber {G}amma(k) in {W}aring's
  {P}roblem}, Proc. London Math. Soc. (2), \textbf{28} (1928), pp.~518--542.

\bibitem{HB92}
{\sc D.~R. Heath-Brown}, {\em The {D}ensity of {Z}eros of {F}orms for which
  {W}eak {A}pproximation {F}ails}, Math. Comp., \textbf{59} (1992),
  pp.~613--623.

\bibitem{HB}
{\sc D.~R. Heath-Brown}, {\em Searching for solutions of {$x^3+y^3+z^3=k$}}, in
  S\'{e}minaire de {T}h\'{e}orie des {N}ombres, {P}aris, 1989--90, vol.~102 of
  Progr. Math., Birkh\"{a}user Boston, Boston, MA, 1992, pp.~71--76.

\bibitem{HBetal}
{\sc D.~R. Heath-Brown, W.~M. Lioen, and H.~J. J.~T. Riele}, {\em On solving
  the {D}iophantine equation {$x^{3} + y^{3} + z^{3} = k$} on a vector
  computer}, Math. Comp., \textbf{61} (1993), pp.~235--244.

\bibitem{IsikKaraOzman20}
{\sc E.~Isik, Y.~Kara, and E.~Ozman}, {\em On ternary {D}iophantine equations
  of signature $(p,p,2)$ over number fields}, Turkish J. Math., \textbf{44}
  (2020), pp.~1197--1211.

\bibitem{isik_kara_ozman_2022}
\leavevmode\vrule height 2pt depth -1.6pt width 23pt, {\em On ternary
  {D}iophantine equations of signature $(p,p,\text{3})$ over number fields}, to
  appear in Canad. J. Math.,  (2022).

\bibitem{JacPintWal}
{\sc M.~J. Jacobson, Jr., A.~Pint\'{e}r, and P.~G. Walsh}, {\em A computational
  approach for solving {$y^2=1^k+2^k+\dots+x^k$}}, Math. Comp., \textbf{72}
  (2003), pp.~2099--2110.

\bibitem{Jarvis08}
{\sc F.~Jarvis and J.~Manoharmayum}, {\em On the modularity of supersingular
  elliptic curves over certain totally real number fields}, J. Number Theory,
  \textbf{128} (2008), pp.~589--618.

\bibitem{JarvisMeekin}
{\sc F.~Jarvis and P.~Meekin}, {\em The {F}ermat equation over
  {$\mathbb{Q}(\sqrt{2})$}}, J. Number Theory, \textbf{109} (2004),
  pp.~182--196.

\bibitem{Kano}
{\sc H.~Kano}, {\em On the equation {$s(1^k+2^k+\cdots +x^k)+r=by^z$}}, Tokyo
  J. Math., \textbf{13} (1990), pp.~441--448.

\bibitem{KaraOzman20}
{\sc Y.~Kara and E.~Ozman}, {\em Asymptotic generalized {F}ermat's last theorem
  over number fields}, Int. J. Number Theory, \textbf{16} (2020), pp.~907--924.

\bibitem{Kellner}
{\sc B.~C. Kellner}, {\em On stronger conjectures that imply the {E}rd{\H
  o}s-{M}oser conjecture}, J. Number Theory, \textbf{131} (2011),
  pp.~1054--1061.

\bibitem{Kempner}
{\sc A.~Kempner}, {\em Bemerkungen zum {W}aringschen {P}roblem}, Math. Ann.,
  \textbf{72} (1912), pp.~387--399.

\bibitem{Kezuka22}
{\sc Y.~Kezuka}, {\em Which numbers are sums of two cubes? (in japanese)}, in
  Proceedings of the Waseda Number Theory Workshop 2022, 2022.

\bibitem{KL20}
{\sc Y.~Kezuka and Y.~Li}, {\em A classical family of elliptic curves having
  rank one and the 2-primary part of their {T}ate-{S}hafarevich group
  non-trivial}, Doc. Math., \textbf{25} (2020), pp.~2115--2147.

\bibitem{KJ22}
{\sc M.~Khawaja and F.~Jarvis}, {\em Fermat's {L}ast {T}heorem over
  {$\mathbb{Q}(\sqrt{2},\sqrt{3})$}},
  \href{https://arxiv.org/abs/2210.03744}{arXiv: 2210.03744},  (2022).

\bibitem{koutsianas2018perfect}
{\sc A.~Koutsianas and V.~Patel}, {\em Perfect powers that are sums of squares
  in a three term arithmetic progression}, Int. J. Number Theory, \textbf{14}
  (2018), pp.~2729--2735.

\bibitem{Kraus}
{\sc A.~Kraus}, {\em Majorations effectives pour l'\'{e}quation de {F}ermat
  g\'{e}n\'{e}ralis\'{e}e}, Canad. J. Math., textbf{49} (1997), pp.~1139--1161.

\bibitem{Kraus19}
\leavevmode\vrule height 2pt depth -1.6pt width 23pt, {\em Le th\'{e}or\`eme de
  {F}ermat sur certains corps de nombres totalement r\'{e}els}, Algebra Number
  Theory, \textbf{13} (2019), pp.~301--332.

\bibitem{kriz20}
{\sc D.~Kriz}, {\em Supersingular main conjectures, {S}ylvester's conjecture
  and {G}oldfeld's conjecture}, preprint,  (2020).

\bibitem{Krzysztofek}
{\sc B.~Krzysztofek}, {\em The equation {$1\sp{n}+2\sp{n}+ \cdots
  +m\sp{n}=(m+1)\sp{n}\cdot k$} (in polish)}, Wy\.{z}. Szko\l . Ped. w
  Katowicach---Zeszyty Nauk. Sekc. Mat., \textbf{5} (1966), pp.~47--54.

\bibitem{kundu2021perfect}
{\sc D.~Kundu and V.~Patel}, {\em Perfect powers that are sums of squares of an
  arithmetic progression}, Rocky Mountain J. Math, \textbf{51} (2021),
  pp.~933--949.

\bibitem{Landau1911}
{\sc E.~Landau}, {\em \"{U}ber die {Z}erlegung positiver ganzer {Z}ahlen in
  positive {K}uben}, Arch. der Math. u. Phys., \textbf{3} (1911), pp.~248--252.

\bibitem{LeSoy2022}
{\sc M.~Le and G.~Soydan}, {\em On the power values of the sum of three squares
  in arithmetic progression}, Math. Commun., \textbf{27} (2022), pp.~137--150.

\bibitem{LeVeque}
{\sc W.~J. LeVeque}, {\em On the equation {$y^{m}=f(x)$}}, Acta Arith.,
  \textbf{9} (1964), pp.~209--219.

\bibitem{Linnik43}
{\sc U.~V. Linnik}, {\em On the representation of large numbers as sums of
  seven cubes}, Rec. Math. [Mat. Sbornik] N. S.,  (1943), pp.~218--224.

\bibitem{Ljunggren}
{\sc W.~Ljunggren}, {\em New solution of a problem proposed by {E}. {L}ucas},
  Norsk Mat. Tidsskr., \textbf{34} (1952), pp.~65--72.

\bibitem{Lucas1180}
{\sc E.~Lucas}, {\em Problem {$1180$}.}, Nouvelle Ann. Math., \textbf{14}
  (1875), p.~336.

\bibitem{LucasSol1180}
\leavevmode\vrule height 2pt depth -1.6pt width 23pt, {\em Solution to
  {P}roblem {$1180$}.}, Nouvelle Ann. Math., \textbf{16} (1877), pp.~429--432.

\bibitem{Lucas1961}
\leavevmode\vrule height 2pt depth -1.6pt width 23pt, {\em Recherches sur
  l'analyse ind\'{e}termin\'{e}e et l'{A}rithm\'{e}tique de {D}iophante},
  Librairie Scientifique et Technique Albert Blanchard, Paris, 1961.

\bibitem{Ma1985}
{\sc D.~G. Ma}, {\em An elementary proof of the solutions to the {D}iophantine
  equation {$6y^2=x(x+1)(2x+1)$}}, Sichuan Daxue Xuebao,  (1985), pp.~107--116.

\bibitem{MaDiophantine}
\leavevmode\vrule height 2pt depth -1.6pt width 23pt, {\em On the {D}iophantine
  equation {$6Y^2=X(X+1)(2X+1)$}}, Kexue Tongbao (English Ed.), \textbf{30}
  (1985), p.~1266.

\bibitem{SonMac1012}
{\sc K.~MacMillan and J.~Sondow}, {\em Divisibility of power sums and the
  generalized {E}rd{\H o}s-{M}oser equation}, Elem. Math., \textbf{67} (2012),
  pp.~182--186.

\bibitem{Mazur78}
{\sc B.~Mazur}, {\em Rational isogenies of prime degree}, Invent. Math,
  \textbf{44} (1978), pp.~129--162.

\bibitem{McCurley84}
{\sc K.~S. McCurley}, {\em An effective seven cube theorem}, J. Number Theory,
  \textbf{19} (1984), pp.~176--183.

\bibitem{michaudjacobs2021fermats}
{\sc P.~Michaud-Jacobs}, {\em Fermat's {L}ast {T}heorem and modular curves over
  real quadratic fields}, Acta Arith., \textbf{203} (2022), pp.~319--351.

\bibitem{Mignotte96}
{\sc M.~Mignotte}, {\em A note on the equation {$ax^n - by^n = c$}}, Acta
  Arith., \textbf{75} (1996), pp.~287--295.

\bibitem{VoutierMignotte}
{\sc M.~Mignotte and P.~M. Voutier}, {\em A kit for linear forms in three
  logarithms (with an appendix by {M}. {L}aurent)},
  \href{https://arxiv.org/abs/2205.08899}{arXiv: 2205.08899},  (2022).

\bibitem{MW54}
{\sc J.~C.~P. Miller and M.~F.~C. Woollett}, {\em Solutions of the
  {D}iophantine equation: {$x^3+y^3+z^3=k$}}, J. Lond. Math. Soc.,
  \textbf{s1-30} (1955), pp.~101--110.

\bibitem{mocanu221}
{\sc D.~Mocanu}, {\em Asymptotic {F}ermat for signatures $(p,p,2)$ and
  $(p,p,3)$ over totally real fields}, Mathematika, \textbf{68} (2022),
  pp.~1233--1257.

\bibitem{mocanu222}
\leavevmode\vrule height 2pt depth -1.6pt width 23pt, {\em Asymptotic {F}ermat
  for signatures $(r,r,p)$ using the modular approach},
  \href{https://arxiv.org/abs/2212.10627}{arXiv: 2212.10627},  (2022).

\bibitem{Mordell56}
{\sc L.~J. Mordell}, {\em On the {I}nteger {S}olutions of the {E}quation
  {$x^2+y^2+z^2+2xyz = n$}}, J. Lond. Math. Soc., \textbf{s1-28} (1953),
  pp.~500--510.

\bibitem{MoreeGeneralised}
{\sc P.~Moree}, {\em Diophantine equations of {E}rd{\H o}s-{M}oser type}, Bull.
  Austral. Math. Soc., \textbf{53} (1996), pp.~281--292.

\bibitem{MoreeRabbits}
\leavevmode\vrule height 2pt depth -1.6pt width 23pt, {\em A top hat for
  {M}oser's four mathemagical rabbits}, Amer. Math. Monthly, \textbf{118}
  (2011), pp.~364--370.

\bibitem{MoreeSurvey}
\leavevmode\vrule height 2pt depth -1.6pt width 23pt, {\em Moser's mathemagical
  work on the equation {$1^k+2^k+\cdots+(m-1)^k=m^k$}}, Rocky Mountain J.
  Math., \textbf{43} (2013), pp.~1707--1737.

\bibitem{MRU1994}
{\sc P.~Moree, H.~J.~J. te~Riele, and J.~Urbanowicz}, {\em Divisibility
  properties of integers {$x,\ k$} satisfying {$1^k+\cdots+(x-1)^k=x^k$}},
  Math. Comp., 63 (1994), pp.~799--815.

\bibitem{MoretBlanc}
{\sc Moret-Blanc}, {\em [solution to] {Q}uestion {$1180$}}, Nouvelle Ann.
  Math., \textbf{15} (1876), pp.~46--48.

\bibitem{Moser1953}
{\sc L.~Moser}, {\em On the diophantine equation $1^n+2^n+3^n+\cdots
  +(m-1)^n=m^n.$}, Scripta Math., \textbf{19} (1953), pp.~84--88.

\bibitem{MullerChabauty}
{\sc J.~S. M\"{u}ller}, {\em Rational points on curves}, 2019.
\newblock
  \hyperlink{https://heilbronn.ac.uk/wp-content/uploads/2019/06/Steffen-Muller-Lecture-Notes-and-Exercises-1-3-updated.pdf}{Online
  Lecture Notes}.

\bibitem{Patel17}
{\sc V.~Patel}, {\em Perfect powers that are sums of consecutive like powers},
  PhD thesis, University of Warwick, 2017.

\bibitem{van18sq}
{\sc V.~Patel}, {\em Perfect powers that are sums of consecutive squares},
  Comptes Rendus Math., \textbf{40} (2018), pp.~33--38.

\bibitem{VP21}
{\sc V.~Patel}, {\em A {L}ucas–{L}ehmer approach to generalised
  {L}ebesgue–{R}amanujan–{N}agell equations}, Ramanujan J., \textbf{56}
  (2021), pp.~585--–596.

\bibitem{patel2017powers}
{\sc V.~Patel and S.~Siksek}, {\em On powers that are sums of consecutive like
  powers}, Res. Number Theory, \textbf{3} (2017), pp.~1--7.

\bibitem{Pillai}
{\sc S.~S. Pillai}, {\em On {W}aring's problem {$g (6)=73$}}, Proc. Indian
  Acad. Sci., Sect. A., \textbf{12} (1940), pp.~30--40.

\bibitem{Pink}
{\sc I.~Pink and Z.~R\'{a}bai}, {\em On the {D}iophantine equation
  {$x^2+5^{k}17^l=y^n$}}, Commun. Math. C, \textbf{19} (2011).

\bibitem{Pinter1997}
{\sc A.~Pint\'{e}r}, {\em A note on the equation
  {$1^k+2^k+\cdots+(x-1)^k=y^m$}}, Indag. Math. (N.S.), \textbf{8} (1997),
  pp.~119--123.

\bibitem{pinter2007power}
{\sc {\'A}.~Pint{\'e}r}, {\em On the power values of power sums}, J. Number
  Theory, \textbf{125} (2007), pp.~412--423.

\bibitem{Csaba}
{\sc {\relax Cs}.~Rakaczki}, {\em On some generalizations of the {D}iophantine
  equation {$s(1^k+2^k+\dots+x^k)+r=dy^n$}}, Acta Arith., \textbf{151} (2012),
  pp.~201--216.

\bibitem{R2007}
{\sc O.~Ramar\'{e}}, {\em An explicit result of the sum of seven cubes},
  Manuscripta Math., \textbf{124} (2007), pp.~59--75.

\bibitem{Ribet90}
{\sc K.~Ribet}, {\em On modular representations of
  {$Gal(\bar{\mathbb{Q}}/\mathbb{Q})$} arising from modular forms}, Invent.
  Math., \textbf{100} (1990), pp.~431--476.

\bibitem{RVZ95}
{\sc F.~Rodr\'{\i}guez~Villegas and D.~Zagier}, {\em Which primes are sums of
  two cubes?}, in Number theory ({H}alifax, {NS}, 1994), vol.~\textbf{15} of
  CMS Conf. Proc., Amer. Math. Soc., Providence, RI, 1995, pp.~295--306.

\bibitem{satge86}
{\sc P.~Satge}, {\em Groupes de {S}elmer et corps cubiques}, J. Number Theory,
  \textbf{23} (1986), pp.~294--317.

\bibitem{Schaffer}
{\sc J.~J. Sch\"{a}ffer}, {\em The equation {$1^p+2^p+3^p+\cdots+n^p=m^q$}},
  Acta Math., \textbf{95} (1956), pp.~155--189.

\bibitem{ShoreyTijdeman86}
{\sc T.~N. Shorey and R.~Tijdeman}, {\em Exponential {D}iophantine equations},
  vol.~\textbf{87} of Cambridge Tracts in Mathematics, Cambridge University
  Press, Cambridge, 1986.

\bibitem{SiegelKum}
{\sc C.~L. Siegel}, {\em Zu zwei {B}emerkungen {K}ummers}, Nachr. Akad. Wiss.
  G{\"o}ttingen Math.-Phys. Kl. II,  (1964), pp.~51--57.

\bibitem{Siegel}
\leavevmode\vrule height 2pt depth -1.6pt width 23pt, {\em \"{U}ber einige
  {A}nwendungen diophantischer {A}pproximationen}, in On some applications of
  {D}iophantine approximations, vol.~\textbf{2} of Quad./Monogr., Ed. Norm.,
  Pisa, 2014, pp.~81--138.

\bibitem{Siksek}
{\sc S.~Siksek}, {\em The {M}odular {A}pproach to {D}iophantine {E}quations},
  vol.~\textbf{36} of Panor. Synth\'{e}ses, Soc. Math, France, Paris, 2012,
  pp.~151--179.

\bibitem{Siksek2015Chabauty}
\leavevmode\vrule height 2pt depth -1.6pt width 23pt, {\em Chabauty and the
  {M}ordell-{W}eil sieve}, vol.~\textbf{41} of NATO Sci. Peace Secur. Ser. D
  Inf. Commun. Secur., IOS, Amsterdam, 2015, pp.~194--224.

\bibitem{Siksek7Cubes}
\leavevmode\vrule height 2pt depth -1.6pt width 23pt, {\em Every integer
  greater than 454 is the sum of at most seven positive cubes}, Algebra Number
  Theory, \textbf{10} (2016), pp.~2093--2119.

\bibitem{SiksekSurvey}
{\sc S.~Siksek}, {\em Sums of integer cubes}, Proc. Natl. Acad. Sci. USA,
  \textbf{118} (2021), p.~e2103697118.

\bibitem{SonMac}
{\sc J.~Sondow and K.~MacMillan}, {\em Reducing the {E}rd{\H o}s-{M}oser
  equation {$1^n+2^n+\dots+k^n=(k+1)^n$} modulo {$k$} and {$k^2$}}, Integers,
  \textbf{11} (2011), pp.~A34, 8.

\bibitem{soydan2017diophantine}
{\sc G.~Soydan}, {\em On the {D}iophantine equation {$(x+1)^k+(x+2)^k+\dots
  +(lx)^k=y^n$}}, Publ. Math. Debrecen, \textbf{91} (2017), pp.~369--382.

\bibitem{Sier60}
{\sc W.~Spierpi\'{n}ski}, {\em On some unsolved problems of arithmetics},
  Manuscripta Math., \textbf{25} (1960), pp.~125--136.

\bibitem{Stroeker1995}
{\sc R.~J. Stroeker}, {\em On the sum of consecutive cubes being a perfect
  square}, Compos. Math., \textbf{97} (1995), pp.~295--307.

\bibitem{TW95}
{\sc R.~Taylor and A.~Wiles}, {\em Ring-theoretic properties of certain {H}ecke
  algebras}, Ann. Math., \textbf{141} (1995), pp.~553--572.

\bibitem{Thue09}
{\sc A.~Thue}, {\em \"{U}ber {A}nn\"{a}herungswerte algebraischer {Z}ahlen.},
  J. Reine Angew. Math., \textbf{135} (1909), pp.~284--305.

\bibitem{Uchiyama}
{\sc S.~Uchiyama}, {\em On a {D}iophantine equation}, Proc. Japan Acad. Ser. A
  Math. Sci., \textbf{55} (1979), pp.~367--369.

\bibitem{Urbanowicz1988}
{\sc J.~Urbanowicz}, {\em Remarks on the equation
  {$1^k+2^k+\cdots+(x-1)^k=x^k$}}, Nederl. Akad. Wetensch. Indag. Math.,
  \textbf{50} (1988), pp.~343--348.

\bibitem{Urbanowicz1994}
\leavevmode\vrule height 2pt depth -1.6pt width 23pt, {\em On {D}iophantine
  equations involving sums of powers with quadratic characters as coefficients.
  {I}}, Compos. Math., \textbf{92} (1994), pp.~249--271.

\bibitem{Urbanowicz1996}
\leavevmode\vrule height 2pt depth -1.6pt width 23pt, {\em On {D}iophantine
  equations involving sums of powers with quadratic characters as coefficients.
  {II}}, Compos. Math., \textbf{102} (1996), pp.~125--140.

\bibitem{van2021sum}
{\sc J.~M. van Langen}, {\em On the sum of fourth powers in arithmetic
  progression}, Int. J. Number Theory, \textbf{17} (2021), pp.~191--221.

\bibitem{VW2000}
{\sc R.~C. Vaughan and T.~D. Wooley}, {\em Waring's problem: a survey}, in
  Number theory for the millennium, {III} ({U}rbana, {IL}, 2000), A K Peters,
  Natick, MA, 2002, pp.~301--340.

\bibitem{VooGyoTij}
{\sc M.~Voorhoeve, K.~Gy{\H o}ry, and R.~Tijdeman}, {\em On the {D}iophantine
  equation {$1^{k}+2^{k}+\cdots +x^{k}+R(x)=y^{z}$}}, Acta Math., \textbf{143}
  (1979), pp.~1--8.

\bibitem{Waldschmidt20}
{\sc M.~Waldschmidt}, {\em Thue {D}iophantine equations -- a survey}, Springer,
  Singapore, (2020), pp.~25--41.

\bibitem{Watson51}
{\sc G.~L. Watson}, {\em A proof of the seven cube theorem}, J. London Math.
  Soc., \textbf{26} (1951), pp.~153--156.

\bibitem{Watson}
{\sc G.~N. Watson}, {\em The problem of the square pyramid.}, Messenger of
  Math., \textbf{48} (1918), pp.~1--22.

\bibitem{Wieferich}
{\sc A.~Wieferich}, {\em Beweis des {S}atzes, da\ss sich eine jede ganze {Z}ahl
  als {S}umme von h\"{o}chstens neun positiven {K}uben darstellen l\"{a}\ss t},
  Math. Ann., \textbf{66} (1908), pp.~95--101.

\bibitem{Wiles}
{\sc A.~Wiles}, {\em Modular elliptic curves and {F}ermat's {L}ast {T}heorem},
  Ann. of Math., \textbf{142} (1995), pp.~443--551.

\bibitem{zhang2014diophantine}
{\sc Z.~Zhang}, {\em On the {D}iophantine equation {$(x- 1)^k+ x^k+(x+ 1)^k=
  y^n$}}, Publ. Math. Debrecen, \textbf{85} (2014), p.~93–100.

\bibitem{zhang2017diophantine4}
\leavevmode\vrule height 2pt depth -1.6pt width 23pt, {\em On the {D}iophantine
  equation {$(x- d)^{4}+ x^4+(x+d)^{4}= y^{n}$}}, Int. J. Number Theory,
  \textbf{13} (2017), pp.~2229--2243.

\bibitem{zhang2017diophantine}
\leavevmode\vrule height 2pt depth -1.6pt width 23pt, {\em On the {D}iophantine
  equations {$(x-1)^{3}+ x^{5}+(x+1)^{3}= y^{n}$} and {$(x- 1)^{5}+
  x^{3}+(x+1)^{5}= y^{n}$}}, Publ. Math. Debrecen, \textbf{91} (2017),
  pp.~383--390.

\bibitem{ZhangBai2013}
{\sc Z.~Zhang and M.~Bai}, {\em On the {D}iophantine equation
  {$(x+1)^2+(x+2)^2+\dots+(x+d)^2=y^n$}}, Funct. Approx. Comment. Math.,
  \textbf{49} (2013), pp.~73--77.

\end{thebibliography}

\end{document}